\documentclass[11pt]{article}

\usepackage{amsbsy,amsfonts,amsmath,amssymb,amsthm
,bbm,enumerate,euscript,graphicx,color}
\usepackage{authblk}
\usepackage{mathtools}
\usepackage[utf8]{inputenc}
\def\R{\mathbb{R}}
\def\N{\mathbb{N}}
\def\Z{\mathbb{Z}}
\def\la{\lambda}
\def\vp{\varphi}

\def\R{\mathbb{R}}
\def\E{\mathbf{E}}
\def\V{\mathbf{Var}}
\def\P{\mathbf{P}}
\def\C{\textbf{Cum}}

\def\rho{\varrho}
\def\D{\mathrm{d}}
\def\I{\mathbf{1}}
\def\vol{|\, \Xi\, \cap\, \rho\, K\, |_2}
\def\uw{w^{\cup}}
\def\vw{w^{\cap}}

\newtheorem{theo}{Theorem}
\newtheorem{coro}{Corollary}

\newtheorem{exam}{Example}
\newtheorem{lemm}{Lemma}
\newtheorem{defn}{Definition}

\theoremstyle{definition}
\newtheorem*{proo}{Proof}
\numberwithin{equation}{section}
\newtheorem*{rema}{Remark}

\title{\bf {On the Variance of the Area of Planar Cylinder Processes
Driven by Brillinger-Mixing Point Processes}}

\author{\textbf{ Daniela Flimmel}\thanks{\noindent Charles University of Prague, daniela.flimmel@karlin.mff.cuni.cz}\,\,\textbf{ and
Lothar Heinrich}{\footnote{University of Augsburg, lothar.heinrich@math.uni-augsburg.de}} }

\date{\today}

\begin{document}

\maketitle

\begin{abstract}
\noindent We study some asymptotic properties of cylinder processes in the plane defined as union sets of dilated straight lines (appearing as mutually overlapping infinitely long strips) derived from a stationary independently marked point process on the real line, where the marks describe thickness and orientation of individual cylinders. Such cylinder processes form an important class of (in general non-stationary) planar random sets. We observe the cylinder process in an unboundedly growing domain $\rho K$ when $\rho \to \infty\,$, where the set $K$ is compact and star-shaped w.r.t. the origin ${\bf o}$ being an inner point of $K$. Provided the unmarked point process satisfies a Brillinger-type mixing condition and the thickness of the typical cylinder has a finite second moment we prove a (weak) law of large numbers as well as a formula of the asymptotic variance for the area of the cylinder process in $\rho K$. Due to the long-range dependencies of the cylinder process, this variance increases proportionally to $\rho^3$.
\end{abstract}

\section{Introduction and Preliminaries}

Cylinder processes in $\R^d$ defined as coutable union of dilated affine subspaces $\R^k\,,\,k=1,\ldots,d-1\,$, are basic random set models in stochastic geometry, see e.g. \cite{Mat75}, \cite{Weil87}, \cite{SchWei08}, \cite{Mol17}. They have numerous applications (mostly for $d = 2, 3$) among others in material sciences  to model materials consisting of long thick fibres, see e.g.  \cite{SpiSpo11}. Until now, so far as we know,  asymptotic properties of cylinder processes in expanding domains were exclusively studied under  Poisson assumptions, see \cite{HeiSpi09},\cite{HeiSpi13}. In this paper, our focus is put on  planar cylinder processes which are derived from stationary independently marked point processes on $\R^1$. Under comparatively strong conditions on the higher-order cumulant measures of the unmarked (ground) point process we are able to prove first, a mean-square limit of the relative part of the area of an expanding star-shaped window covered by the union of cylinders, and second, we derive an explicit formula for the asymptotic variance of this area. The latter is an important first step in proving asymptotic normality of the covered area that will be carried out in a later paper. Our main results, Theorems 1 and 2, in Section 2 generalize some of the results obtained in \cite{HeiSpi13} (in particular Theorem 2 in \cite{HeiSpi13}) for stationary Poisson cylinder processes even under general dimensional assumptions.

\smallskip\noindent
Throughout in this paper, all random elements are defined on a common probability space $[\Omega, \mathcal{F}, \P]$ and by $\E$ resp. $\V$, we denote the expectation resp. variance w.r.t. $\P$. Next we describe a cylinder process in $\R^2$ in terms of its generating stationary, independently marked point process on $\R^1$. 
For doing this, let $(\Phi_0, R_0)$ be the generic random vector taking value in the mark space $[0,\pi]\times[0,\infty)$ that describes the orientation $\Phi_0$ and the cross-section (or base) $\Xi_0:=[-R_0,R_0]$ of the typical cylinder. In addition, we assume that $R_0 \sim F$ and $\Phi_0 \sim G$ are independent, i.e. $\P(R_0\le r,\Phi_0\le \vp) = F(r)\,G(\vp)$. Now we introduce a stationary independently marked point process as locally finite, simple counting measure $\Psi_{F,G}^P := \sum_{i \in \Z}\delta_{[P_i,(\Phi_i,R_i)]}$ defined on the Borel sets of $\R^1\times[0,\pi]\times[0,\infty)\,$, whose finite-dimensional distributions are shift-invariant in the first component,  see e.g. \cite{DVJ03}, \cite{Fraetal81} or \cite{SchWei08}. The stationary unmarked (or ground) point process $\Psi=\sum_{i \in \Z}\delta_{P_i} \sim P$ with finite and positive intensity $\la=\E\Psi([0,1])>0$ is assumed to be independent of the i.i.d. sequence $\{(\Phi_i,R_i): i \in \Z:=\{0,\pm 1, \pm 2,\ldots\}\}$ of mark vectors. Each triplet $[P_i,(\Phi_i,R_i)],i \in \Z,$ determines a random cylinder $g(P_i,\Phi_i)\oplus b(\mathbf{o},R_i)\,$, 
where $b(\mathbf{o},r)$ is the circle in $\R^2$ with radius $r \geq 0$ and centre in the origin $\mathbf{o}$ and $\oplus$ stands for pointwise addition (Minkowski sum) of subsets of $\R^2$. Here, $g(p,\vp):=\{(x,y) \in \R^2: x\,\cos\vp + y\,\sin\vp=p\}$ denotes the unique line with signed distance $p \in \R^1$ from $\mathbf{o}$ and an angle $\vp \in [0,\pi)$  measured anti-clockwise between the normal vector $v(\vp)=(\cos\vp,\sin\vp)^T$ on the line with direction in the half plane not containing $\mathbf{o}$ and the $x$-axis.

\smallskip\noindent
The intensity measure  $\Lambda_{F,G}((\cdot)\times[0,\vp]\times[0,r]) := \E\Psi_{F,G}^P((\cdot)\times[0,\vp]\times[0,r])$ of $\Psi_{F,G}^P$ can be expressed  for any $\;r \ge 0$ and $0 \le \vp \le \pi$ as
$$ \Lambda_{F,G}((\cdot)\times[0,\vp]\times[0,r]) = \E\Psi(\cdot)\,\P(\Phi_0\le \vp ,R_0 \le r) = \la\,|\cdot|_1\;G(\vp)\,F(r)\quad\mbox{with}\quad \la > 0\,,$$
where $|\cdot|_k$ denotes the Lebesgue measure on $\R^k$. 
Now we are in a position to define the main subject of this paper. 

\begin{defn}
A \textit{cylinder process} $\Xi = \Xi_{F,G}^P$ in the Euclidean plane $\R^2$ derived from the stationary independently marked point process $\Psi_{F,G}^P$  is defined by random union set
\begin{equation}\label{CylinderProces}
\Xi_{F,G}^P := \bigcup_{i \in \Z} \big(\,g(P_i,\Phi_i)\oplus b(\mathbf{o},R_i)\,\big)\,,
\end{equation}
which in general is neither closed nor stationary.
\end{defn}

\smallskip\noindent
For more details and  a general survey on cylinder processes we refer to \cite{Weil87}, see also the monographs \cite{Mat75}, \cite{Mol17}. 
The aim of this paper consists first, in proving the $L^2-$convergence of the ratio $|\Xi \cap \rho K|_2/|\rho K|_2$ to a deterministic limit as $\rho \to \infty$ and second, in proving the existence and determining the explicit shape of the asymptotic variance 
\begin{equation}\label{eq1.2}
\lim\limits_{\rho \to \infty}\frac{\V(|\Xi \cap \rho K|_2)}{\rho^3}=: \sigma^2_P(K,F,G)\,,
\end{equation}
for some fixed compact star-shaped set $K \subset \R^2$ containing the origin ${\bf o}$ as inner point. The limit $\sigma^2_P(K,F,G)$ is positive and finite (if $\E|\Xi_0|^2 = 4\,\E R_0^2 < \infty $) and  depends on the shape of $K$, the first and second moment of $F$ and the distribution function $G$ which is assumed to be continuous (not necessarily absolutely continuous). A purely discrete distribution function $G$ yields different expressions for $\sigma^2_P(K,F,G)$ even if  $\Psi \sim P =\Pi_{\la}$ is a stationary Poisson point process with intensity $\lambda > 0$, see \cite{HeiSpi09},\cite{HeiSpi13}. A distribution function $G$ without jumps implies that $\P(\Phi_0=\Phi_1) = 0$ if the angles $\Phi_0, \Phi_1 \sim G$ are independent. 

\smallskip\noindent
Note that the order $\rho^3$ of the growth of $\V(|\Xi \cap \rho K|_2)$ is much faster that the growth of the area $|\rho K|_2=\rho^2|K|_2$ which reveals a typical feature of long-range dependencies within the random set \eqref{CylinderProces}.

\smallskip\noindent
We recall the fact that the probability space $[\Omega, \mathcal{F}, \P]$ on which the marked point process $\Psi_m$ is defined can be chosen in such a way that the mapping $(x,\omega)\mapsto \mathbf{1}_{\Xi(\omega)}(x) \in \{0,1\}$ for $(x,\omega) \in \R^2\times\Omega$ is measurable w.r.t. the product $\sigma$-field $\mathcal{B}(\R^2)\otimes\mathcal{F}$, see Appendix in \cite{Hei05}. This enables us to apply Fubini’s theorem to the random field of indicator variables $\{\mathbf{1}_{\Xi}(x), x \in \R^2\}$ and implies that the $k$th-order {\it mixed moment function} 

\begin{equation}\label{MixedMoments}
p^{(k)}_{\Xi}(x_1,\dots,x_k):=\E\Big(\prod_{j=1}^k \mathbf{1}_{\Xi}(x_j)\Big)=\P(x_1\in \Xi,\dots,x_k\in\Xi),\quad x_1,\dots,x_k\in \R^2,
\end{equation}
are $\mathcal{B}(\R^{2k})-$measurable for any $k \in \N :=\{1,2,\ldots\}\,$.

\smallskip\noindent
The distribution of a random closed set $\Xi$ is determined by its {\it Choquet functional}
\begin{equation}\label{Choquet}
T_{\Xi}(X) := \P(\Xi \cap X \neq \emptyset)\quad\mbox{for}\quad X \in \mathcal{K}_2,
\end{equation}
where  $\mathcal{K}_2$ denotes the family of non-empty compact sets in $\R^2$.
In particular, the $k$th order moment functions $p_{\Xi^c}^{(k)}$ of the $0-1$-random field $\xi(x):=\mathbf{1}_{\Xi^c}(x)$ can be expressed by \eqref{MixedMoments} and \eqref{Choquet} for any $k\geq 1$:
$$p^{(k)}_{\Xi^c}(x_1,\ldots,x_k)= \E\Big(\prod_{j=1}^k \xi(x_j)\Big)=\P(\{x_1,\ldots,x_k\}\cap \Xi =\emptyset)=1-T_{\Xi}(\{x_1,\ldots,x_k\}).$$
\begin{lemm}\label{Lemma1}
For any $X \in \mathcal{K}_2$, we have
\begin{equation}\label{eq1.5}
 T_{\Xi}(X) = 1 - G_P\big[1-\P\big((\cdot) \in [-R_0,R_0]\oplus \langle v(\Phi_0),X \rangle\big)\big],   
\end{equation}
where $\langle v(\Phi_0),X \rangle:= \bigcup_{x \in X} \langle v(\Phi_0),x \rangle$ with scalar product $\langle \cdot, \cdot \rangle$ in $\R^2$ and $G_P[w(\cdot)]$ denotes the \textit{probability generating functional} (short: pgf) of $\Psi \sim P$ defined for Borel-measurable functions $w: \R^1 \to [0,1]$ by
\begin{equation}\label{pgf}
G_P[w(\cdot)] := \E\Big( \prod_{i: \Psi(\{P_i\}) > 0} w(P_i)\Big), \quad \text{where } \int_{\R^1}(1-w(x)) \D x < \infty.
\end{equation}
\end{lemm}

\begin{coro}\label{Coro1}
For $X =\{x_1,\ldots,x_k\}$ with pairwise distinct points $x_1\ldots, x_k \in \R^2$ we get 
\begin{equation}\label{eq1.7}
T_{\Xi^c}(\{x_1,\ldots,x_k\})=G_P\Big[1-\P\Big((\cdot)\in \bigcup_{i=1}^k([-R_0,R_0] + \langle v(\Phi_0), x_i\rangle)\Big)\Big].
\end{equation}
\end{coro}

\begin{exam}\label{Exam1}
For a stationary Poisson process $\Psi \sim \Pi_{\la}$ with intensity $\la > 0$, we have $G_{\Pi_{\la}}[w(\cdot)]=\exp\{\la \int_{\R^1}(w(x)-1)\D x\}$ implying that
\begin{align*}
    T_{\Xi}(X)& = 1-\exp\Big\{ -\la \int_{\R^1}\P\left((g(p,\Phi_0)\oplus b(\mathbf{o},R_0))\cap X \neq \emptyset\right) \D p \Big\}\\
    &= 1- \exp\Big\{-\la \E \big|[-R_0,R_0]\oplus \langle v(\Phi_0),X\rangle\big|_1\Big\}\\
    &= 1 -\exp\Big\lbrace  -\la \int_0 ^{\infty}\int_{0}^{\pi}\big|[-r,r]\oplus \langle v(\vp),X\rangle\big|_1 \D G(\vp) \D F(r)\Big\rbrace.
\end{align*}
\end{exam}

\noindent 
In the special case $X =\{x_1,\ldots,x_k\}$ such that $x_i=(x_i^{(1)}, x_i^{(2)})^T$ and $\langle v(\vp),x_i \rangle = x_i^{(1)} \cos\,\vp + x_i^{(2)} \sin\,\vp$ for $i=1,\ldots,k\,$, it follows from 
Corollary \ref{Coro1} that
\begin{align*}
    T_{\Xi}(\{x_1,&\ldots, x_k\}) = 1 - \exp\Big\lbrace -\la \E \big| \bigcup_{i=1}^k ([-R_0,R_0]+x_i^{(1)}\cos(\Phi_0 + x_i^{(2)}\sin\,\Phi_0) \big|_1\,\Big\rbrace\\
    & = 1 - \exp\Big\lbrace -\la \int_0^{\infty}\int_{0}^{\pi}\big| \bigcup_{i=1}^k([-r,r]+ x_i^{(1)} \cos\,\vp +x_i^{(2)}\sin\,\vp)\big|_1\D G(\vp)\D F(r) \Big\rbrace.
\end{align*}

\bigskip\noindent
\begin{proo}[Lemma \ref{Lemma1}]
To prove formula \eqref{eq1.2}, we need the orthogonal matrix
\begin{equation}\label{orthmatrix}
O(\varphi) = \begin{pmatrix}
\cos\,\vp & -\sin\,\vp\\
\sin\,\vp & \cos\,\vp
\end{pmatrix},
\end{equation}
which represents an anti-clockwise rotation by the angle $\varphi\in[0,\pi)$ so that $O(-\varphi)v(\varphi)=(1,0)^T$ and $O(\varphi)(1,0)^T = v(\varphi)$ since it holds $O(-\varphi) = O^T(\varphi)=O^{-1}(\varphi)\,$.
Using the pgf \eqref{pgf} and the independence assumption in the definition of \eqref{CylinderProces}, we obtain
\begin{align}\label{eq1.9}
T_{\Xi}(X)  &= 1 -\P\big(\Xi \cap X = \emptyset\big)
    = 1 - \P\Big(\bigcap_{i: \Psi(\{P_i\}) > 0}\{(g(P_i,\Phi_i)\oplus b(\bf{o},R_i)) \cap X = \emptyset\}\Big)\nonumber \\
    & =1 - \E \Big(\prod_{i:\Psi(\{P_i\}) > 0} \mathbf{1}_{\{(g(P_i,\Phi_i)\oplus b(\bf{o},R_i)) \cap X = \emptyset\}}\Big) \nonumber \\
    &= 1- \int\limits_{\mathbf{N}}\E \Big(\prod_{i: \psi(\{p_i\}) > 0}\big(1 - \P\big((g(p_i,\Phi_i)\oplus b(\mathbf{o},R_i)) \cap X = \emptyset \big|\Psi = \psi \big)\big)\Big)\P(\Psi \in \D\psi)\nonumber \\
    &= 1- \int\limits_{\mathbf{N}} \prod_{i: \psi(\{p_i\}) > 0}\big(1 - \P\big((g(p_i,\Phi_i)\oplus b(\mathbf{o},R_i)) \cap X = \emptyset\big)\big)\P(\Psi \in \D\psi)\nonumber \\
    &=1- \int\limits_{\mathbf{N}} \prod_{i: \psi(\{p_i\}) > 0}\big(1 - \P\big(p_i \in [-R_0,R_0]\oplus \langle v(\Phi_0),X\rangle\big)\big)\P(\Psi \in \D\psi),
\end{align}
where $\mathbf{N}$ denotes the set of locally finite simple counting measures   on the Borel-$\sigma-$algebra $\mathcal{B}(\R^1)$.  The last step leading to \eqref{eq1.9} is seen as follows:
\begin{align*}
    \Big\{\big(\,g(p,\Phi_0)\oplus b(\mathbf{o},R_0)\,\big)\cap X \neq \emptyset \Big\} &=\Big\{ p\, v(\Phi_0) \in \big(-g(0,\Phi_0)\oplus b(\mathbf{o},R_0)\big)\oplus X \Big\}\\
    &=\Big\{ p\, O(-\Phi_0)\,v(\Phi_0) \in \big(g(0,0)\oplus b(\mathbf{o},R_0)\big)\oplus O(-\Phi_0)X \Big\}\\
    &=\Big\{ p\, (1,0)^T \in \big(g(0,0)\oplus b(\mathbf{o},R_0)\big)\oplus O(-\Phi_0)X \Big\}\\
    &=\Big\{p \in [-R_0, R_0] \oplus \langle v(\Phi_0), X\rangle \Big\}.
\end{align*}
Obviously, \eqref{eq1.9} coincides with \eqref{eq1.5}. Hence, the proof of Lemma 1 is complete. \qed
\end{proo}

\section{Factorial moment expansion of $\E |\Xi \cap \rho K|_2$ and $\V |\Xi \cap \rho K|_2$}
The proof of our asymptotic results relies on an  expansion of the pgf \eqref{pgf} (resp. its logarithm) in terms of the factorial moment (resp. cumulant) measures, see Chapter 5.5 in \cite{DVJ03} or \cite{BlaMerSch97}. To begin with, let us fix $K \subset \R^2$ to be a compact, star-shaped set containing the origin $\mathbf{o}$ as an inner point. Further, let $\rho \geq 1$ be a scaling factor tending to infinity implying that $\rho K \uparrow \R^2$ as $\rho \to \infty$. The second-order mixed moment functions  \eqref{MixedMoments} fulfill the relation
\begin{equation}\label{eq2.1}
    p_{\Xi}^{(2)}(x_1,x_2)-p_{\Xi}^{(1)}(x_1)\, p_{\Xi}^{(1)}(x_2) = p_{\Xi^c}^{(2)}(x_1,x_2) - p_{\Xi^c}^{(1)}(x_1)\,p_{\Xi^c}^{(1)}(x_2).
\end{equation} 
By applying Fubini's theorem we get, together with \eqref{eq1.5},  that
\begin{align}\label{eq2.2}
    \E|\Xi \cap \rho K|_2 &=\E \int\limits_{\R^2} \I_{\Xi}(x) \,\I_{\rho K}(x)\D x = \int\limits_{\rho K}p_{\Xi}^{(1)}(x)\D x = \rho^2 \int\limits_K T_{\Xi}(\{\rho  x\})\D x\nonumber\\
    &=\rho^2\int\limits_K  \big(\,1- G_P\big[1-\P\big((\cdot)\in [-R_0,R_0] + \rho \langle v(\Phi_0), x\rangle\,\big)\big]\big)\D x.
\end{align}
For the variance, we get from \eqref{MixedMoments}, \eqref{eq1.7} and \eqref{eq2.1}  that
$$\V\big(|\Xi \cap \rho K|_2\big)=\int\limits_{\rho K}\int\limits_{\rho K}\big(p_{\Xi^c}^{(2)}(x_1,x_2)-p_{\Xi^c}^{(1)}(x_1)\, p_{\Xi^c}^{(1)}(x_2)\big)\D x_1 \D x_2.$$
Together with \eqref{eq1.7}, we obtain to the following lemma.
\begin{lemm}\label{Lemma2}
With the above notation and $\langle v(\vp),x_i\rangle = x^{(1)}_i \cos \vp + x^{(2)}_i \sin \vp$ for $i = 1,2$, we have
\begin{align}\label{eq2.3}
    \V\big(|\Xi \cap \rho K|_2\big) = \int\limits_{\rho K}\int\limits_{\rho K} &\Big(G_P\big[1-\P\big((\cdot)\in \bigcup_{i=1}^2([-R_0,R_0] + \langle v(\Phi_0), x_i\rangle)\big)\big]\nonumber\\
    &- \prod_{i=1}^2 G_P\big[1-\P\big( (\cdot)\in [-R_0,R_0] + \langle v(\Phi_0), x_i\rangle\big)\big]\Big)\D x_1 \D x_2.
\end{align}
\end{lemm}

\medskip\noindent
Formula \eqref{eq2.3} can be generalized to higher-order cumulants $\C_k\big(|\Xi \cap \rho K|_2\big)$  for any $k \geq 3$, where the $k$th-order cumulant $\C_k(X)$ of a random variable $X$ can be expressed by its moments 
$\E X,\ldots,\E X^k$ as follows:
  
\begin{equation}
\C_k(X) = \sum_{\ell=1}^k (-1)^{\ell-1}(\ell-1)!\sum\limits_{\substack{K_1\cup\cdots\cup K_\ell\\=\{1,\dots,k\}}}\prod_{j=1}^\ell \E X^{\#K_j}=k!\sum_{\ell=1}^k\frac{(-1)^{\ell-1}}{\ell}\sum\limits_{\substack{k_1+\cdots + k_\ell=k \\ k_j \geq 1,j=1,\ldots,\ell}}\prod_{j=1}^\ell \frac{\E X^{k_j}}{k_j!}\,,
\nonumber\end{equation}
where the first inner sum runs over all decompositions of $\{1,\dots,k\}$ into $\ell$ disjoint non-empty subsets $K_1,\dots,K_\ell$ and $\# K_i$ denotes the number of elements of $K_i, i=1,\dots,\ell.$

\medskip\noindent 
Combining the latter representation  and the  formula

$$\E|\Xi \cap \rho K|^k_2 = \E \int\limits_{(\rho K)^k}\prod_{i=1}^k \I_{\Xi}(x_i)\D(x_1,\ldots,x_k)=\int\limits_{(\rho K)^k} p^{(k)}_{\Xi}(x_1, \ldots, x_k)\D(x_1,\ldots,x_k)\;\;\mbox{for}\;\;k \geq 2,$$

\smallskip\noindent  
with the $k$th-order {\it mixed cumulant function} of the random field $\{\I_{\Xi}(x),x\in\R^2\}$
$$
c_{\Xi}^{(k)}(x_1,\ldots,x_k) := \sum_{\ell=1}^k (-1)^{\ell-1}(\ell-1)! \sum
\limits_{\substack{K_1\cup\cdots\cup K_\ell\\=\{1,\dots,k\}}}
\prod_{j=1}^{\ell}\,p_{\Xi}^{(\# K_j)}(x_i:i\in K_\ell)\;\;\mbox{for}\;\; k \in \N\,,
\nonumber$$
which satisfy the identity $c_{\Xi}^{(k)}(x_1,\ldots,x_k) = (-1)^k\;c_{\Xi^c}^{(k)}(x_1,\ldots,x_k)$ for $k\ge 2$, we arrive at 
\begin{equation}\label{eq2.4}
\qquad \C_k(|\Xi \cap \rho\,K|_2) = (-1)^k\;\int\limits_{(\rho\,K)^k}\;c_{\Xi^c}^{(k)}(x_1,\ldots,x_k)\,{\rm d}(x_1,\ldots,x_k)\;\;\mbox{for}\;\;k\ge 2\,,  
\end{equation}
where 
$$ p_{\Xi^c}^{(k)}(x_1,\ldots,x_k)= 1-T_{\Xi}(\{x_1,\ldots,x_k\})=G_P\bigl[1-{\bf P}\bigl((\cdot)\in \bigcup_{i=1}^k([-R_0,R_0] + \langle v(\Phi_0),x_i\rangle)\bigr)\bigr]\,.\nonumber
$$


\noindent
In order to treat the moments and cumulants of $\vol$, the following relations are useful. Let $a_1, a_2,\ldots$ be real numbers in $[0,1]$. Then we have
\begin{equation}\label{eq2.5}
    1-\prod_{i=1}^n(1-a_i)=\sum_{k=1}^n\frac{(-1)^{k-1}}{k!}\sum\limits_{1\leq i_1,\ldots,i_k\leq n}^{\neq} a_{i_1}\cdot\ldots\cdot a_{i_k}\quad \text{for } n\geq 1.
\end{equation}
Moreover, for any odd number $m<n$ (provided $n\geq 2$), the so called \textit{Bonferroni inequalities} (see e.g. \cite{GalSim96}) hold:
\begin{equation}\label{eq2.6}
\sum_{k=1}^{m+1}\frac{(-1)^{k-1}}{k!}\sum\limits_{1\leq i_1,\ldots,i_k\leq n}^{\neq} a_{i_1}\cdot\ldots\cdot a_{i_k}\leq 1-\prod_{i=1}^n (1-a_i)\leq \sum_{k=1}^{m}\frac{(-1)^{k-1}}{k!}\sum\limits_{1\leq i_1,\ldots,i_k\leq n}^{\neq} a_{i_1}\cdot\ldots\cdot a_{i_k}.
\end{equation}

\begin{defn}\label{notation}
To simplify the notation, we define  for $k \ge 2$ (not necessarily pairwise distinct) points $x_1,\ldots,x_k \in \R^2$ and $\Xi_0 = [-R_0,R_0]$ the functions
\begin{equation*}
\uw_{x_1,\ldots,x_k}(p):= \P\Big(p\in \bigcup_{i=1}^k(\Xi_0 + \langle v(\Phi_0), x_i\rangle)\Big)\;\;\;\mbox{and}\;\;\;
\vw_{x_1,\ldots,x_k}(p):= \P\Big(p\in \bigcap_{i=1}^k(\Xi_0 + 
\langle v(\Phi_0), x_i \rangle)\Big)\,.
\end{equation*}
For $k=1$ we put $\uw_{x}(p)=\vw_{x}(p):=w_{x}(p)$. Obviously, $\uw_{x_1,x_2}(p)=w_{x_1}(p)+w_{x_2}(p)-\vw_{x_1,x_2}(p)\,$.
\end{defn}

\medskip\noindent 
As a consequence of  \eqref{eq2.2} and \eqref{eq2.6} and the definition of the factorial moment measures $\alpha^{(k)}(\cdot)$ of $\Psi\sim P$, we get the following series expansion
\begin{align}\label{eq2.7}
    \E\vol = & \int_{\rho\,K} \big(1-G_P\big[1- w_x(\cdot)\big]\big) \D x = \rho^2 \,\int_K \big(1-G_P\big[1- w_{\rho x}(\cdot)\big]\big) \D x\nonumber\\
    = &\rho^2 \sum_{k=1}^{\infty}\frac{(-1)^{k-1}}{k!}\int_K\int_{\R^k}\prod_{j=1}^k w_{\rho x}(p_j)\,\alpha^{(k)}(\D(p_1,\ldots,p_k))\D x,
\end{align}
provided that the infinite sum on the right hand side converges. From \eqref{eq2.6} we obtain immediately the estimates
\begin{align}\label{eq2.8}
\bigg|\;1-G_P[1-w_{\rho x}(\cdot)] -&\sum_{k=1}^{m-1}\frac{(-1)^{k-1}}{k!}\int_{\R^k} \prod_{j=1}^k w_{\rho x}(p_j) \alpha^{(k)}(\D(p_1,\ldots,p_k))\; \bigg|\nonumber \\
    \leq &\;\;\frac{1}{m!}\int_{\R^m}\prod_{j=1}^m w_{\rho x}(p_j)\,\alpha^{(m)} (\D(p_1,\ldots,p_m))
\end{align}
for any $m \geq 1$. It is easily seen that the right hand side of \eqref{eq2.7} is convergent if and only if
\begin{equation}\label{eq2.9}
\frac{1}{m!} \int_{\R^m}\,\prod_{j=1}^m w_{\rho x}(p_j)\,\alpha^{(m)} (\D(p_1,\ldots,p_m)) \xrightarrow[m \to \infty]{} 0\,.
\end{equation}

\noindent
One way to show \eqref{eq2.9} consists of expressing $\alpha^{(m)}(\cdot)$ by \textit{factorial cumulant measures} $\gamma^{(k)}(\cdot), k=1,\ldots,m$ where $\gamma^{(1)}(B)=\alpha^{(1)}(B)=\la\,|B|_1$ and for $k \geq 2$,
\begin{equation}\label{CumulantMeasureReprez}
\alpha^{(k)}(\times_{i=1}^k B_i)=\sum_{\ell=1}^k \sum\limits_{\substack{K_1\cup\cdots\cup K_\ell\\=\{1,\dots,k\}}}\prod_{j=1}^\ell \gamma^{(\#K_j)}(\times_{i \in K_j}B_i).
\end{equation}

\noindent
The representation \eqref{CumulantMeasureReprez} follows by inverting the defining formula for $\gamma^{(k)}(\cdot)$ which is as follows:
$$ \gamma^{(k)}(\times_{i=1}^k B_i):= \sum_{\ell=1}^k (-1)^{\ell-1}(\ell-1)!\sum\limits_{\substack{K_1\cup\cdots\cup K_\ell\\=\{1,\dots,k\}}}\prod_{j=1}^\ell \alpha^{(\#K_j)}(\times_{i \in K_j}B_i)\,.$$
The latter formula is based on the general relationship  between mixed moments and mixed cumulants, see \cite{LeoShi59} or \cite{SauSta91}. Note that $\gamma^{(k)}$ is a locally finite, signed measure on $[\R^k, \mathcal{B}(\R^k)]$.

\medskip\noindent
Due to the stationarity of $\Psi \sim P$ we can implicitely define the $k$th- order \textit{reduced cumulant measure} $\gamma^{(k)}_{red}(\cdot)$ as the unique signed measure on $[\R^{k-1},\mathcal{B}(\R^{k-1})]$ satisfying
$$ \gamma^{(k)}(\times_{i = 1}^k B_i)= \la\,\int_{B_1} \gamma^{(k)}_{red}(\times_{i=2}^k(B_i-p))\D p$$
for all bounded sets $B_1,\ldots,B_k \in \mathcal{B}(\R^1)$. 

\medskip\noindent
The \textit{total variation measure} $|\gamma_{red}^{(k)}|(\cdot)$ is defined by $|\gamma_{red}^{(k)}|(\cdot)=(\gamma_{red}^{(k)})^+(\cdot)+(\gamma_{red}^{(k)})^-(\cdot)$, where the measures $(\gamma_{red}^{(k)})^+(\cdot)$ and $(\gamma_{red}^{(k)})^-(\cdot)$ are given by the Jordan decomposition of the signed measure $\gamma_{red}^{(k)}(\cdot)=(\gamma_{red}^{(k)})^+(\cdot)-(\gamma_{red}^{(k)})^-(\cdot)$. The \textit{total variation} of $\gamma_{red}^{(k)}(\cdot)$ on $[\R^{k-1}, \mathcal{B}(\R^{k-1})]$ is defined by 
$\|\gamma_{red}^{(k)}\|_{TV}:= |\gamma_{red}^{(k)}|(\R^{k-1})$. Furthermore, if $\gamma_{red}^{(k)}(\cdot)$ possesses a Lebesgue density $c_{red}^{(k)}(\cdot)$ on $\R^{k-1}$ (called $k$th-order \textit{reduced cumulant density}), we define the canonical $L_q$-norm $\|c_{red}^{(k)}\|_q := \big(\,\int_{\R^{k-1}}|c_{red}^{(k)}(x)|^q {\rm d}x\,\big)^{1/q}$ for $k\ge 2$ and the modified $L_q^*$-norm  $\|c_{red}^{(k)}\|_q^* := \int_{\R^1}\big(\int_{\R^{k-2}}|c_{red}^{(k)}(x,p)|^q {\rm d}x\big)^{1/q}{\rm d}p$ for $k\ge 3$, where $1 \le q < \infty\,$. Formally we may put $\|\gamma_{red}^{(1)}\|_{TV}:= 1$ and $\|c_{red}^{(1)}\|_q:= 1$.

\begin{defn}\label{Brillinger}
A stationary point process $\Psi \sim P$ on $[\R^1,\mathcal{B}(\R^1)]$ with intensity $\la =\E\Psi([0,1]) > 0$ satisfying $\E \Psi^{k}([0,1])<\infty$ for all $k\geq 2$, is called Brillinger-mixing  if $\|\gamma_{red}^{(k)}\|_{TV} < \infty$ for all $k \geq 2\,$.
$\Psi \sim P$ is said to be strongly Brillinger-mixing $($ strongly $L_q-$Brillinger-mixing, resp. strongly $L_q^*-$Brillinger-mixing for some $q \ge 1\,)$ if there are constants $b > 0$ and $a \geq b^{-1}$ such that $\|\gamma^{(k)}_{red}\|_{TV}\leq a\, b^k \,k!$    
$($ if $c_{red}^{(k)}(\cdot)$ exists such that $\|c_{red}^{(2)}\|_1 < \infty$  and $\|c^{(k)}_{red}\|_q \le a_q\, (b_q)^k \,k!$ for $k \ge 2$ with constants $a_q, b_q > 0$   resp. $\|c^{(k)}_{red}\|_q^* \le a_q^*\, (b_q^*) ^k \,k!$  for $k \ge 3$ with constants $a_q^*, b_q^* > 0$\,$)$.
\end{defn}

\smallskip\noindent
\begin{rema}\label{Rem2}
For formal reason we put $\|\gamma^{(1)}_{red}\|_{TV}:=1$ so that $a \geq b^{-1}$ makes $\|\gamma^{(1)}_{red}\|_{TV}:=1 \le a\,b\,$. Further, note that the existence and integrability of $c_{red}^{(k)}(\cdot)$ imply that $\|c_{red}^{(2)}\|_1 = \|\gamma_{red}^{(2)}\|_{TV}$ and $\|c_{red}^{(k)}\|_1 = \|c_{red}^{(k)}\|^*_1 = \|\gamma_{red}^{(k)}\|_{TV}$ for all $k \geq 3\,$. 
\end{rema}

\medskip\noindent
\begin{rema}\label{Rem3}
In general, the Brillinger-mixing condition is formulated for stationary point processes on $\R^d, d \ge 1$. This condition expresses some kind of weak correlatedness (or asymptotic uncorrelatedness) of the numbers of points lying in bounded sets having a large (or unboundedly increasing) distance of each other. This type of weak dependence does not necessarily imply ergodicity, see \cite{Hei18}, but allows to prove central limit theorems for various stochastic models related with point processes, e.g. in stochastic geometry, statistical physics for $d\ge 1$ or in queueing theory for $d=1$, see e.g. \cite{HeiSch85}. In \cite{HeiPaw13, Hei21} the relations between (strong) Brillinger-mixing and classical mixing conditions are studied. Strong Brillinger-mixing requires exponential moments of the number of points in  bounded sets. For any dimension $d \ge 1$, examples of such point processes are determinental point processes, see \cite{Hei16}, \cite{BisLav16}, Poisson cluster processes if the number of daughter points has an exponential moment and certain Cox processes as well as Gibbsian processes under suitable restrictions, see \cite{Rue69}. For $d=1$, renewal processes with an exponentially decaying interrenewal density, see \cite{HeiSch85}, among them the Erlang process and the Macchi process, see \cite{DVJ03} (p. 144), are strongly Brillinger-mixing.   
\end{rema}

\bigskip\noindent
\begin{lemm}\label{Lemma3}
If the stationary point process $\Psi \sim P$ is strongly Brillinger-mixing with $b < \frac{1}{2}$ and $\E R_0 < \infty\,$, then 
\begin{equation}\label{eq2.11}
  \sum_{m=1}^\infty  \frac{1}{m!} \,\int_{\R^m}\prod_{j=1}^m w_{\rho x}(p_j) \alpha^{(m)}(\D(p_1,\ldots,p_m))\, \le \frac{b}{1-2b}\,\big(\exp\{a\,\la\,\E|\Xi_0|_1\}-1\big)\,,
\end{equation}
which immediately implies \eqref{eq2.9}. If $\Psi \sim P$ is strongly $L_q-$Brillinger-mixing for some $q > 1$ such that $b_q < \frac{1}{2}\,(\E |\Xi_0|)^{\frac{1}{q}-1}\,$, then the estimate \eqref{eq2.11} remains valid with $a$ and $b$ replaced by $a_q\,(\E|\Xi_0|)^{\frac{1}{q}-1}$ and $b_q\,(\E|\Xi_0|)^{1-\frac{1}{q}}\,$, respectively.
\end{lemm}

\bigskip\noindent
\begin{proof}
Using the representation \eqref{CumulantMeasureReprez}, we obtain
\begin{align*}
    &\frac{1}{m!} \int\limits_{\R^m}\prod_{j=1}^m w_{\rho x}(p_j) \alpha^{(m)} (\D(p_1,\ldots,p_m))\\
    &= \frac{1}{m!} \sum_{\ell=1}^m \sum\limits_{\substack{K_1\cup\cdots\cup K_\ell\\=\{1,\dots,m\}}}\prod_{j=1}^\ell \int\limits_{\R^{\# K_j}}
    \prod_{i\in K_j} w_{\rho x}(p_i) \gamma^{(\# K_j)}(\D(p_i: i\in K_j))
\end{align*}

\begin{equation}\label{eq2.12}
=\frac{1}{m!} \sum_{\ell=1}^m \frac{1}{\ell!}\sum\limits_{\substack{k_1+\cdots+k_\ell=m\\k_i\geq 1, i=1\dots,\ell}}\frac{m!}{k_1!\cdots k_\ell!}\prod_{j=1}^\ell f(k_j) = \sum_{\ell=1}^m \frac{1}{\ell!}\sum\limits_{\substack{k_1+\cdots+k_\ell=m\\k_i \geq 1, i=1\dots,\ell}}\prod_{i=1}^\ell\frac{f(k_i)}{k_i!}\,,
\end{equation}   
where $$f(k) := \int\limits_{\R^k}\prod_{i=1}^k w_{\rho x}(p_i) \gamma^{(k)}(\D(p_1,\dots,p_k)) = \la \int\limits_{\R^1}w_{\rho x}(p_1)\int\limits_{\R^{k-1}}\prod_{i=2}^k w_{\rho x}(p_i+p_1)\gamma^{(k)}_{red}(\D(p_2,\ldots,p_k))\D p_1$$ for $k=1,\ldots,m$.
The equality \eqref{eq2.12} is justified by the invariance of  $\gamma^{(k)}(\times_{i=1}^kB_i)$ against permutations of the bounded sets $B_1,\ldots,B_k \in \mathcal{B}(\R^1)$ for each $k \in \N$. We proceed with
\begin{align*}
    |f(k)| \leq \,\la \int_{\R^1}w_{\rho x}(p_1)\int_{\R^{k-1}}\big|\gamma^{(k)}_{red}\big|(\D(p_2,\ldots,p_k))\,\D p_1 = \la\,\E|\Xi_0|_1\,\| \gamma^{(k)}_{red} \|_{TV} \le a\,\la\,\E|\Xi_0|_1\,b^k\,k!\,.
\end{align*}
Here, we have used Fubini's theorem combined with $w_{\rho x}(p)\leq 1$ for 
$p \in \R^1$ and $x \in \R^2$ so that
$$\int_{\R^1}w_{\rho x}(p)\D p=\int_{\R^1}\P\big(p \in \Xi_0 +\rho \langle v(\Phi_0),x\rangle\big)\D p = \int_{\R^1}\P\big(p \in \Xi_0\big)\D p = \E |\Xi_0|_1\,.$$

\smallskip\noindent
Hence, together with the combinatorial relations 
\begin{equation*}\sum\limits_{\substack{k_1+\cdots+k_\ell=m\\ k_i\geq 1, i=1\dots,\ell}} 1 ={{m-1}\choose{\ell-1}}\quad\mbox{and}\quad \sum_{\ell=1}^m  {{m-1}\choose{\ell-1}}=2^{m-1}
\end{equation*}
we arrive at
\begin{align}\label{eq2.13}
\sum_{\ell=1}^m \frac{1}{\ell!}\,\sum\limits_{\substack{k_1+\cdots+k_\ell=m\\k_i\geq 1, i=1\dots,\ell}}\prod_{i=1}^\ell\frac{|\,f(k_i)\,|}{k_i!} & \leq b^m \sum_{\ell=1}^m \frac{(a\,\la \E|\Xi_0|_1)^\ell}{\ell!}{{m-1}\choose{\ell-1}} \leq b^m \;2^{m-1} \; \max_{1\leq \ell 
    \leq m} \frac{(a\,\la \E|\Xi_0|_1)^\ell}{\ell!}\nonumber\\
 & \leq \frac{1}{2}\,\big(\exp\{a\,\la\,\E|\Xi_0|_1\}-1\big)\;(2\,b)^m\,.
\end{align}
By combining \eqref{eq2.12} and \eqref{eq2.13} with $b < 1/2$ the relation \eqref{eq2.11} follows immediately. Under the strong $L_q-$Brillinger-mixing condition we may rewrite $f(k)$ for  $k \ge 2$ as follows:
$$f(k) = \la \int\limits_{\R^1}w_{\rho x}(p_1)\,\E\int\limits_{\R^{k-1}}\prod_{i=2}^k \I_{\Xi_i+\rho\langle v(\Phi_i),x\rangle}(p_i+p_1)\,c^{(k)}_{red}(p_2,\ldots,p_k)\,\D(p_2,\ldots,p_k)\D p_1\,,$$ 
where $\Xi_i=[-R_i,R_i]$ and $(R_2,\Phi_2),\ldots,(R_k,\Phi_k)$  are i.i.d. random vectors with same distribution as $(R_0.\Phi_0)$. Applying H\"older's inequality for $q > 1$ and $p=q/(q-1)\,$, Lyapunov's inequality $\E|\Xi_0|^\frac{1}{p} \le (\E|\Xi_0|)^\frac{1}{p} = (\E|\Xi_0|)^{1-\frac{1}{q}}$  and the condition $\|c^{(k)}_{red}\|_q \le a_q\,b_q^k\,k!$, we obtain that

\begin{align*}
|f(k)| \le \la \|c^{(k)}_{red}\|_q\, \E |\Xi_1|_1\;\prod_{i=2}^k\E|\Xi_i|^{\frac{1}{p}}_1  \leq \la \|c^{(k)}_{red}\|_q\;(\E|\Xi_0|_1)^{1+\frac{k-1}{p}} \le \la\,a_q\,(\E|\Xi_0|_1)^\frac{1}{q}\;(b_p\,(\E|\Xi_0|_1)^{1-\frac{1}{q}})^k\;k!\,.
\end{align*}

\smallskip\noindent
By repeating the foregoing steps with the latter bound the proof of Lemma 3 is finished. \end{proof}

\begin{lemm}\label{Lemma4}
Let $\Psi\sim P$ be a stationary point process on $\R^1$ satisfying $max_{2\le k \le m}\|\gamma_{red}^{(k)}\|_{TV} < \infty$ for some fixed $m \ge 2\,$. If  $\,\E R_0 < \infty$ and $\Phi_0 \sim G$ has a continuous distribution function $G\,$ then, for $m \ge 2$ not necessarily distinct point $x_1,\ldots,x_m \in \R^2 \setminus \{\bf{o}\}\,$, 
\begin{equation}\label{eq2.14}
\int\limits_{\R^m}\prod_{j=1}^m w_{\rho x_j}(p_j)\;\alpha^{(m)}(\D(p_1,\ldots,p_m)) \xrightarrow[\rho \to \infty]{} \la^m\;\prod_{j=1}^m \int\limits_{\R^1}w_{\rho x_j}(p)\,\D p = (\la\,\E|\Xi_0|_1)^m\,. 
\end{equation}
\end{lemm}
\medskip
\begin{proo}[Lemma \ref{Lemma4}] We use the representation \eqref{CumulantMeasureReprez} for $k=m$ to rewrite the difference of left-hand and right-hand side of  \eqref{eq2.14} as follows:

\medskip
\begin{equation*}
\sum_{\ell=1}^{m-1} \sum\limits_{\substack{K_1\cup\cdots\cup K_\ell\\=\{1,\dots,m\}}}\;\prod_{j=1}^\ell\; \int\limits_{\R^{\#K_j}} \prod_{i\in K_j} w_{\rho x_i}(p_i)\gamma^{(\#K_j)}({\rm d}(p_i:i \in K_j))\,.
\end{equation*}
Hence, the limit \eqref{eq2.14} is shown if and only if the finite  sum in the latter line disappears as $\rho \to \infty$ and this in turn follows by showing that, for $k=2,\ldots,m\,$,
$$
\int\limits_{\R^k} \prod_{i=1}^k w_{\rho x_i}(p_i)\gamma^{(k)}({\rm d}(p_1,\ldots,p_k)) = \la \int\limits_{\R^k}w_{\rho x_1}(p_1)\prod_{i=2}^k w_{\rho x_i}(p_i+p_1)\,\gamma_{red}^{(k)}({\rm d}(p_2,\ldots,p_k)){\rm d}p_1 \xrightarrow[\rho \to \infty]{} 0\,.
$$
In view of $0 \le w_{\rho x_i}(p_i+p_1) \le 1$ for $i=3,\ldots,k$ it is sufficient to prove 
$$
\int\limits_{\R^{k-1}}\int\limits_{\R^1}\P(p_1\in \Xi_0+\rho \langle v(\Phi_0),x_1 \rangle)\,\P(p_1\in \Xi_0+\rho \langle v(\Phi_0),x_2 \rangle -p_2)\,{\rm d}p_1\,|\gamma_{red}^{(k)}|({\rm d}(p_2,\ldots,p_k))\xrightarrow[\rho \to \infty]{} 0\,.
$$
Since the total variation measure $|\gamma_{red}^{(k)}|(\cdot)$ is bounded on $\R^{k-1}$ and the inner integral over $\R^1$ is less than or equal to $\E|\Xi_0|_1\,$, we have only to verify that the inner integral disappears as $\rho \to \infty\,$. For this purpose, we rewrite its integrand as expectation
$\E\,\I_{\{\Xi_1+\rho \langle v(\Phi_1),x_1 \rangle\}}(p_1)\,\I_{\{\Xi_2+\rho \langle v(\Phi_2),x_2 \rangle -p_2\}}(p_1)\,$,
where $\Xi_i=[-R_i,R_i]$ and $\Phi_i$ for $i=1,2$ have the same distribution as  $\Xi_0=[-R_0,R_0]$ and $\Phi_0$, respectively, and $R_1,R_2,\Phi_1,\Phi_2$ are independent of each other. By Fubini's theorem and the shift-invariance of the Lebesgue measure, we arrive at
\begin{align*}
 & \int\limits_{\R^1}\P(\,p_1\in \Xi_0+\rho \langle v(\Phi_0),x_1 \rangle\,)\;\P(\,p_1\in \Xi_0+\rho \langle v(\Phi_0),x_2 \rangle -p_2\,)\,{\rm d}p_1\nonumber \\
= & \; \int\limits_{\R^1}\E\,\I_{\{\Xi_1+\rho \langle v(\Phi_1),x_1 \rangle\}}(p_1)\,\I_{\{\Xi_2+\rho \langle v(\Phi_2),x_2 \rangle - p_2\}}(p_1)\,{\rm d}p_1 \nonumber\\
= & \; \E\big|\,\Xi_1\cap\big(\Xi_2-p_2+\rho ( \langle v(\Phi_2),x_2 \rangle - \langle v(\Phi_1),x_1 \rangle )\big)\,\big|_1\xrightarrow[\rho \to \infty]{} 0\,.
\end{align*}

The limit in the last line can verified as follows: We fix two points $x_i=\|x_i\|(\cos(\alpha_i),\sin(\alpha_i))\in\R^2,\,i=1,2\,$, and two points  $v(\vp_i)=(\cos(\vp_i),\sin(\vp_i)),\,i=1,2\,$, on the unit circle line. It is easily seen that the equation $\langle v(\vp_1),x_1 \rangle = \langle v(\vp_2),x_2 \rangle$, i.e. $\|x_1\|\cos(\vp_1-\alpha_1)=\|x_2\|\cos(\vp_2-\alpha_2)$ holds for at most a finite number of pairs $\vp_1,\vp_2 \in [0,\pi]\,$. Hence, for two independent random angles $\Phi_1, \Phi_2$ with common atomless distribution function $G(\cdot)$ we have 
$$ \P( \langle v(\Phi_1),x_1 \rangle \ne \langle v(\Phi_2),x_2 \rangle ) = 1 \quad\mbox{for any two points}\quad x_1, x_2 \in \R^2\quad\mbox{with}\quad  \|x_1\|+\|x_2\| > 0\,.$$   \qed
\end{proo}

\medskip\noindent
From  Lemma \ref{Lemma4} and \eqref{eq2.8} we obtain the behaviour  of the expectation  of $|\Xi\cap \rho\,K|_2$ as $\rho \to \infty\,$.
\begin{coro}\label{Coro2}
Let $\Psi\sim P$ be a Brillinger-mixing point process $\R^1$. If $\,\E R_0 < \infty$ and $\Phi_0 \sim G$ has a continuous distribution function $G\,$ then  
$$\frac{\E\vol}{|\rho K|_2} \xrightarrow[\rho \to \infty]{}\sum_{k=1}^{\infty}\frac{(-1)^{k-1}}{k!}(\la\,\E|\Xi_0|_1)^k = 1-\exp\{-\la\,\E|\Xi_0|_1\}\,.$$
\end{coro}

\begin{proo}[Corollary \ref{Coro2}] An application of \eqref{eq2.14} for $x_1=\cdots=x_m=x \ne {\bf{o}}$ to the inequality \eqref{eq2.8} yields

\begin{equation}\label{eq2.15} 
\bigg|\;\lim_{\rho\to\infty}\big(1 - G_P[1-w_{\rho x}(\cdot)]\big) -\sum_{k=1}^{m-1}\frac{(-1)^{k-1}\,(\la\,\E|\Xi_0|_1)^k}{k!}\;\bigg| 
    \leq \frac{(\la\,\E|\Xi_0|_1)^m}{m!}\quad\mbox{for any}\quad m \ge 1\,.
\end{equation}
Combining this with \eqref{eq2.7} leads to 
$$\frac{\E\vol}{|\rho K|_2} = \frac{1}{|K|_2}\int\limits_K\big(1 - G_P[1-w_{\rho x}(\cdot)]\big){\rm d}x \xrightarrow[\rho \to \infty]{} \sum_{k=1}^{\infty}\frac{(-1)^{k-1}\,(\la\,\E|\Xi_0|_1)^k}{k!}$$
which immediately gives the assertion of Corollary \ref{Coro2}. \qed
\end{proo}

\section{Main Results}

\bigskip
The first result can be considered as a  planar mean-square ergodic theorem which implies a weak law of large numbers for $\Xi \cap \rho K$ in the Euclidean plane $\R^2$. 

\bigskip
\begin{theo}\label{theorem1}
Assume that the stationary point process $\Psi \sim P$ on $\R^1$ is  Brillinger-mixing. Further suppose that $\E R_0 < \infty$ and $\Phi_0 \sim G$ has a continuous distribution function $G$. Then

\bigskip
\begin{equation}\label{eq3.1}
\E \bigg( \frac{\, \vol\, }{\, |\, \rho\, K\,|_2} - \big(\,1 - \exp\{-\la\,\E\, |\, \Xi_0\,|_1\}\,\big)\bigg)^2 \xrightarrow[\rho \to \infty]{}0\quad\mbox{with}\quad \Xi_0 := [-R_0, R_0]\,.
\end{equation}

\end{theo}

\bigskip\noindent
Our second result provides  the exact asymptotic behavior of the variance of the area of the cylinder process (1.1) that is contained in a star-shaped set $\rho K$ which is growing unboundedly in all directions. For this purpose, in comparison with Theorem 1, we need a  strengthening and quantification of the classical Brillinger-mixing condition. 

\bigskip
\begin{theo}\label{theorem2}
Assume that the stationary point process $\Psi \sim P$ on $\R^1$  is either strongly Brillinger-mixing with $b < 1/2$ or  strongly $L_q-$Brillinger-mixing with $(\E |\Xi_0|_1)^{1-\frac{1}{q}}\;b_q < 1/2$ and strongly $L^*_q-$Brillinger-mixing with $(\E |\Xi_0|_1)^{1-\frac{1}{q}}\;b^*_q < 1/2$ for some $q > 1\,$, where $\Xi_0 := [-R_0, R_0]$. Further suppose that $\E R_0 ^2 < \infty$ and $\Phi_0 \sim G$ has a continuous distribution function $G$. Then

\bigskip
\begin{equation}\label{eq3.2}
\lim\limits_{\rho \to \infty}  \frac{\V(\,\vol\, )}{\rho^3} = \la\, e^{-2\,\la \,\E\, |\, \Xi_0\, |_1}\, \Big(\, (\, \E\, |\,\Xi_0\, |_1\, )^2 \;\gamma_{red}^{(2)}(\R^1)\; C_1^{G,K} + 2\, \E\, |\, \Xi_0\, |_1^2\; C_2^{G,K}\,\Big)\,,
\end{equation}

\medskip\noindent
where 
\begin{equation*}
C_1 ^{G,K} :=\int\limits_{\R^1} (\,\E|\,g(p, \Phi_0) \cap K\,|_1\,)^2\, \D p\quad\mbox{and}\quad
C_2 ^{G,K} := \int\limits_0^{\pi} \int\limits_0 ^{r_K(\vp\pm \pi/2)} |\,K \cap (K + s\, v(\vp \pm \frac{\pi}{2})\,|_2 \,\D s \,\D G(\vp)\,.
\end{equation*}
\end{theo}

\bigskip
\begin{rema}\label{remark 4}
In the special case $K=b({\bf{o}},1)$, one can show that $C_1^{G,K}=\frac{16}{3}$ and  $C_2^{G,K} = \frac{8}{3}$ are independent of the distribution function $G$. If $\Phi_0$ is uniformly distributed on $[0,\pi]$, then we get

\begin{align*}
C_1 ^{G,K}& = \frac{1}{\pi^2}\,\int_{-\infty}^{\infty}\Big(\int_0^\pi  |\,g(p,\vp) \cap K\,|_1 \D \vp\,\Big)^2\,\D p\;,\\ 
C_2 ^{G,K}& = \frac{1}{2\pi}\int\limits_0^{2\pi}\int\limits_0^{\infty}|\,K\cap (K + s\,v(\vp)\,)\,|_2\,\D s\,\D\vp
=\frac{1}{2\pi}\int\limits_{\R^2}|\,K \cap (K+x)\,|_2\,\frac{\D x}{\|x\|} = \frac{1}{2\pi}\int\limits_K \int\limits_K \frac{\D x \D y}{\|x-y\|}.
\end{align*}

\bigskip\noindent
The latter double integral is known as \textit{second-order chord power integral of $K$}, see e.g. \cite{HeiSpi13}, p. 327, and  \cite{SchWei08}, Chapt. 7, for integral geometric background.
\end{rema}

\bigskip\noindent
The following two lemmas are essential for the calculation of the right-hand side of \eqref{eq3.2}. Interestingly, the assumptions to prove these lemmas are rather mild in comparison with the Brillinger-mixing-type conditions  in the Theorems 1 and 2.

\bigskip\noindent 
\begin{lemm}\label{lemma5}
Let $\Psi \sim P$ be a second-order stationary point process on $\R^1$ satisfying $\|\gamma_{red}^{(2)}\|_{TV} < \infty$. Further, suppose that $\E R_0 < \infty$ and $\Phi_0 \sim G$ with a not necessarily continuous distribution function $G$. Then 

\bigskip
\begin{align}\label{eq3.3}
    \rho\int\limits_K\int\limits_K \int\limits_{\R^2} w_{\rho x_1}(p_1)\, w_{\rho x_2}(p_2)& \,\gamma^{(2)}\,(\D (p_1,p_2))\,\D x_1 \,\D x_2 \\
    &\xrightarrow[\rho \to \infty]{}\la\,(\E\,|\,\Xi_0\,|_1)^2\,\gamma_{red}^{(2)}(\R^1)\,\int\limits_{\R^1}\,\big(\,\E\,|\,g(p,\Phi_0)\cap K\,|_1\,\big)^2\,\D p\,.\nonumber
\end{align}
\end{lemm}

\bigskip\noindent 
\begin{lemm}\label{lemma6}
Assume that $\E R_0^2 < \infty$ and $\Phi_0 \sim G$ with  a not necessarily continuous distribution function $G$. Then 
$$\rho \int\limits_K \int\limits_K \int\limits_{\R^1} \vw_{\rho x_1, \rho x_2}(p)\,\D p\, \D x_1\, \D x_2\xrightarrow[\rho \to \infty]{} 2\,\E\,|\,\Xi_0\,|_1^2\int\limits_0^{\pi}\int\limits_0^{r_K(\vp\pm\frac{\pi}{2})}|\,K\cap(K+s\,v(\vp \pm \frac{\pi}{2}))\,|_2\,\D s\,\D G(\vp)$$
with $r_K(\psi):=\max\{r \geq 0: r v(\psi) \in K \oplus (-K)\}$. Obviously, it holds $r_K(\psi) = r_K(\psi \pm \pi)$.
\end{lemm}

\bigskip
\begin{rema}\label{remark5} Note that in Theorem 1 and 2  the interval $\Xi_0 := [-R_0,R_0]$ with $\E R_0^k < \infty$ can be replaced by a finite union of random closed intervals $\Xi_0 \subset \R^1$ satisfying $\inf \Xi_0\leq 0 \leq \sup \Xi_0$  and $\E|\Xi_0|_1^k < \infty$ for $k=1$ or $k=2\,$, respectively. This restriction is based on the definition of a process of
 cylinders with non-convex bases, see e.g. \cite{SpiSpo11}. In Lemma 5 and 6  the cross section (or base) $\Xi_0$ of the typical cylinder can be chosen as random compact set satisfying $0 < \E|\Xi_0|_1 < \infty$ or $\E|\Xi_0|_1^2 < \infty\,$, respectively. 
\end{rema}

\section{Proofs of the Main Results}

\medskip\noindent
\begin{proo}[Theorem \ref{theorem1}] The expectation on the left-hand side of \eqref{eq3.1} can be expressed as follows: 
$$
\frac{\V(|\Xi \cap \rho K|_2)}{|\rho\,K|_2^2} + \bigg(\frac{\E|\Xi \cap \rho K|_2}{|\rho\,K|_2}\,-\,\Big(\,1 - \exp\{- \la\,\E|\Xi_0|_1\}\,\Big)\bigg)^2\,.
$$
In view of Corollary \ref{Coro2} it remains to prove that $\rho^{-4}\,\V(|\Xi \cap \rho K|_2)\xrightarrow[\rho \to \infty]{}0\,$. Using Lemma \ref{Lemma2} and the notation introduced in Section 2 we get

\begin{align*}
\rho^{-4}\,\V(|\Xi \cap \rho K|_2) = \rho^{-4}\,\int\limits_{\rho K}\int\limits_{\rho K}\Big(G_P\big[1-\uw_{x_1,x_2}(\cdot)\big]- \prod_{i=1}^2 G_P\big[1-w_{x_i}(\cdot)\big]\Big)\,{\rm d}x_1{\rm d}x_2\\
= 
\int\limits_K\int\limits_K\Big(G_P\big[1-\uw_{\rho x_1,\rho x_2}(\cdot)\big]\;-\; G_P\big[1-w_{\rho x_1}(\cdot)\big] \;G_P\big[1-w_{\rho x_2}(\cdot)\big]\Big)\,{\rm d}x_1{\rm d}x_2\,.
\end{align*}
Thus, we just have to show that the integrand disappears as $\rho \to \infty$ for distinct points $x_1,x_2\in K\setminus\{{\bf o}\}$, that is, 
\begin{align}\label{eq4.1}
 G_P\big[1-\uw_{\rho x_1,\rho x_2}(\cdot)\big]\;- \;G_P\big[1-w_{\rho x_1}(\cdot)\big]\;G_P\big[1-w_{\rho x_2}(\cdot)\big] \xrightarrow[\rho \to \infty]{} 0\,.
\end{align}
We make use of the finite expansion \eqref{eq2.8} of the pgf $G_P\big[1-w_{\rho x}(\cdot)\big]$ with remainder term, where $w_{\rho x}$ can be replaced by any Borel-measurable function $w : \R^1 \mapsto [0,1]\,$. For brevity, we put 
$$
S_m(w) := \sum_{k=0}^{m-1}(-1)^kT_k(w)\;\;\mbox{with}\;\;T_0(w) := 1\;\;\mbox{and}\;\;T_k(w) := \frac{1}{k!}\int\limits_{\R^k}\prod_{j=1}^kw(p_j)\alpha^{(k)}(d(p_1,\ldots,p_k))$$ 
for $1 \le k \le m \in \N\,$. Hence, \eqref{eq2.8} reads as $\big|\,G_P\big[1-w(\cdot)\big] - S_m(w)\,\big| \le T_m(w)$ which leads us to the following estimate
\begin{align}\label{eq4.2}
&\Big|\,G_P\big[1-\uw_{\rho x_1,\rho x_2}(\cdot)\big]- G_P\big[1-w_{\rho x_1}(\cdot)\big]\; G_P\big[1-w_{\rho x_2}(\cdot)\big] -\Big(S_m(\uw_{\rho x_1,\rho x_2}) - S_m(w_{\rho x_1})\,S_m(w_{\rho x_2})\Big)\,\Big|  
\nonumber \\
&\nonumber\\
&\le \Big|\,G_P\big[1-\uw_{\rho x_1,\rho x_2}(\cdot)\big] - S_m(\uw_{\rho x_1,\rho x_2})\,\Big|+\Big|\,G_P\big[1-w_{\rho x_1}(\cdot)\big] - S_m(w_{\rho x_1})\,\Big|\,G_P\big[1-w_{\rho x_2}(\cdot)\big] \nonumber\\
&\qquad \qquad\qquad\qquad\qquad\qquad\qquad\qquad\;\;\quad + \Big|\,G_P\big[1-w_{\rho x_2}(\cdot)\big] - S_m(w_{\rho x_2})\,\Big| \,S_m(w_{\rho x_1})    \nonumber\\
&\le  T_m(\uw_{\rho x_1,\rho x_2})+T_m(w_{\rho x_1})+T_m(w_{\rho x_2}) + T_m(w_{\rho x_1})\,T_m(w_{\rho x_2})\quad\mbox{for}\quad m \ge 2\,.
\end{align}
Here, we have additionally used that $G_P\big[1-w(\cdot)\big]\le 1$ and   
$\big|S_m(w)\big| \le G_P\big[1-w(\cdot)\big]+T_m(w)\,.$

\noindent
We are now in a position to apply the limit  \eqref{eq2.12} under the assumptions of Lemma \ref{Lemma4}. This yields for $i=1,2$ and $m \in \N$
$$
T_m(w_{\rho x_i})\xrightarrow[\rho \to \infty]{} \frac{(\la\,\E|\Xi_0|_1)^m}{m!}\;\;\mbox{and}\;\;S_m(w_{\rho x_i})\xrightarrow[\rho \to \infty]{} \sum_{k=0}^{m-1}\frac{(-\la\,\E|\Xi_0|_1)^k}{k!}= e^{-\la\,\E|\Xi_0|_1}+\theta_1\,\frac{(\la\,\E|\Xi_0|_1)^m}{m!}
$$
for some $\theta_1\in [-1,1]$ in accordance with $\Big|e^{-x}-\sum_{k=0}^{m-1}\frac{(-x)^k}{k!}\Big|\le \frac{x^m}{m!}$ for any $m \in \N$ and $x \ge 0\,$.

\noindent
Next, we have to find the limit  of $T_m(\uw_{\rho x_1,\rho x_2})$ as $\rho \to \infty\,$. Using the relation $\uw_{x_1,x_2}(p)=w_{x_1}(p)+w_{x_2}(p)-\vw_{x_1,x_2}(p)$ and taking into account that the factorial moment measure $\alpha^{(m)}$ is invariant under permutation of its $m$ components,
we may write 
\begin{align}\label{eq4.3}
T_m(\uw_{\rho x_1,\rho x_2})& = \frac{1}{m!} \int\limits_{\R^m}\prod_{j=1}^m \big(w_{\rho x_1}(p_j)+w_{\rho x_2}(p_j)-\vw_{\rho x_1,\rho x_2}(p_j)\big)\,\alpha^{(m)}(\D(p_1,\ldots,p_m))\nonumber\\
& = \frac{1}{m!} \int\limits_{\R^m}\prod_{j=1}^m \big(w_{\rho x_1}(p_j)+w_{\rho x_2}(p_j)\big)\,\alpha^{(m)}(\D(p_1,\ldots,p_m))\\
& + \frac{1}{m!} \sum_{\ell=1}^m {m \choose \ell}\int\limits_{\R^m}\prod_{i=1}^\ell \vw_{\rho x_1,\rho x_2}(p_i)\prod_{j=\ell+1}^m 
(w_{\rho x_1}(p_j)+w_{\rho x_2}(p_j))\,\alpha^{(m)}(\D(p_1,\ldots,p_m))\nonumber.
\end{align}
There is at least one term $\vw_{\rho x_1,\rho x_2}(p_i)=
\P\big(p_i\in (\Xi_0 +\rho\langle v(\Phi_0),x_1\rangle)\cap (\Xi_0 +\rho\langle v(\Phi_0),x_2\rangle)\big)$ in each summand of  the last line which will be integrated over $\R^1$ w.r.t. ${\rm d}p_i$ so that after expressing $\alpha^{(m)}$ by cumulant measures, see \eqref{CumulantMeasureReprez}, the expectation $\E|\Xi_0 \cap (\Xi_0 +\rho\,\langle v(\Phi_0),x_2-x_1\rangle)|_1$ emerges and disappears as $\rho \to \infty$ if $x_1 \ne x_2$. Thus , the last line disappears completely as  $\rho \to \infty$, whereas the line \eqref{eq4.3}
converges to the limit $(2\,\la\,\E|\Xi_0|_1)^m/m!$ as $\rho \to \infty$ by applying the limit \eqref{eq2.13} once more. 
Therefore, we obtain for any $m \in \N$ that $T_m(\uw_{\rho x_1,\rho x_2})\xrightarrow[\rho \to \infty]{} (2\,\la\,\E|\Xi_0|_1)^m/m!$ and
$$
S_m(\uw_{\rho x_1,\rho x_2})\xrightarrow[\rho \to \infty]{} \sum_{k=0}^{m-1}\frac{(-2\,\la\,\E|\Xi_0|_1)^k}{k!}= e^{-2\,\la\,\E|\Xi_0|_1}+\theta_2\,\frac{(2\,\la\,\E|\Xi_0|_1)^m}{m!}\;\;\mbox{for some}\;\;\theta_2\in [-1,1]. 
$$
The latter limit combined with above limits of $S_m(w_{\rho x_i})$ for $i=1,2$ leads to
$$ 
\overline{\lim_{\rho \to \infty}}\Big|\,S_m(\uw_{\rho x_1,\rho x_2}) - S_m(w_{\rho x_1})\,S_m(w_{\rho x_2})\,\Big|\le \frac{(2\,\la\,\E|\Xi_0|_1)^m}{m!}+2\,\frac{(\la\,\E|\Xi_0|_1)^m}{m!}+\frac{(\la\,\E|\Xi_0|_1)^{2m}}{(m!)^2}\,.
$$
For any given $\varepsilon \in (0,1]$ we find some $m(\varepsilon)$  such that $\frac{(2\,\la\,\E|\Xi_0|_1)^m}{m!} \le \varepsilon$ for all $m \ge m(\varepsilon)\,$.

\medskip\noindent
Thus, the right-hand side of the last inequality does not exceed  $2\,\varepsilon+\varepsilon^2$ for sufficiently large $m\,$. The same bound can be obtained for the limit (as $\rho \to \infty$) of the four summands in line \eqref{eq4.2}. 

\medskip\noindent
Finally, after summarizing all $\epsilon-$bounds of the above limiting terms we arrive at 
$$
\overline{\lim_{\rho \to \infty}}\Big|\,G_P\big[1-\uw_{\rho x_1,\rho x_2}(\cdot)\big]- G_P\big[1-w_{\rho x_1}(\cdot)\big]\; G_P\big[1-w_{\rho x_2}(\cdot)\big]\,\Big| \le 2\,(2\,\varepsilon + \varepsilon^2) \le 6\,\varepsilon\,.
$$

\medskip\noindent
This implies \eqref{eq4.1} completing  the proof of Theorem 1. \qed

\end{proo}

\bigskip\noindent
\begin{proo}[Lemma \ref{lemma5}]
By the stationarity of $\Psi \sim P$ we may write $\gamma^{(2)}({\rm d}(p_1,p_2))=\la\, \gamma^{(2)}_{red}({\rm d}p_2-p_1)\,{\rm d}p_1$ which gives 
\begin{equation*} 
\rho\int\limits_{\R^2}\int\limits_K\int\limits_K w_{\rho x}(p_1)\,w_{\rho y}(p_2)\,\D x \,\D y \,\gamma^{(2)}(\D(p_1,p_2))=
    \rho\,\la\int\limits_{\R^2}\int\limits_K\int\limits_K w_{\rho x}(p_1)\,w_{\rho y}(p_2+p_1)\, \D x\, \D y \,\gamma^{(2)}_{red}(\D p_2)\,\D p_1
\end{equation*}

\noindent
To determine the limit of the right hand side  as $\rho \to \infty$, we rewrite the probabilities $w_{\rho x}(p_1)=\P(p_1\in \{\cdots\})$ and $w_{\rho y}(p_2+p_1) = \P(p_2+p_1\in\{\cdots\})$ by means of the expectation (as integral over the product of probability measures) over the corresponding indicator function $\I_{\{\cdots\}}$. We fix $\Xi_i=\xi_i$ (compact sets in $\R^1$) and $\Phi_i=\varphi_i$ (angles in $[0,\pi]$) for $i=1,2$ and omit the expectation which stands in front of all other integrals due to Fubini's theorem. The intensity $\la$ will be suppressed. Further, we write $x=(x_1, x_2)$ and $y=(y_1, y_2)$. Thus, we only treat the integral

\begin{align}\label{eq4.4}
   & \rho \int\limits_{\R^2}\int\limits_K\int\limits_K\I_{\xi_1+\rho(x_1\cos\vp_1 + x_2\sin\vp_1)}(p_1) \I_{\xi_2+\rho(y_1\cos\vp_2 + y_2\sin\vp_2)}(p_2+p_1)\, \D(x_1,x_2)\,\D(y_1, y_2)\,\gamma^{(2)}_{red}(\D p_2)\,\D p_1\nonumber\\
&= \rho \int\limits_{\R^2}\int\limits_{\R^2}\int\limits_{\R^2}\I_K(x_1, x_2)\I_K(y_1, y_2)\I_{\xi_1+\rho(x_1\cos\vp_1 + x_2\sin\vp_1)}(p_1) \I_{\xi_2+\rho(y_1\cos\vp_2 + y_2\sin\vp_2)}(p_2+p_1)\nonumber \\
   &\hspace{4cm}\times\,\D(x_1,x_2)\,\D(y_1, y_2)\, \gamma^{(2)}_{red}(\D p_2)\,\D p_1 =: J_\rho(K,\xi_1,\vp_1,\xi_2,\vp_2). 
\end{align}

\smallskip\noindent
Now, we substitute $(x_1,x_2)^T = O(\vp_1)(u_1, u_2)^T, (y_1,y_2)^T = O(\vp_2)(v_1, v_2)^T$ , where $O(\vp_1)$ and $O(\vp_2)$ are defined by \eqref{orthmatrix}. Then $x_1= u_1\cos\vp_1 - u_2\sin\vp_1,\; x_2= u_1\sin\vp_1 + u_2\cos\vp_1$ and $y_1= v_1\cos\vp_2 - v_2\sin\vp_2,\; y_2= v_1\sin\vp_2 + v_2\cos\vp_2$. Hence, since $O(\vp_i)^{-1} = O(-\vp_i)$ for $i=1,2$, the integral $ J_\rho(K,\xi_1,\vp_1,\xi_2,\vp_2)$ in \eqref{eq4.4} takes on the form

\begin{align*}
    &\phantom{==}\rho \int\limits_{\R^2}\int\limits_{\R^2}\int\limits_{\R^2}\I_{O(-\vp_1)K}(u_1, u_2)\I_{O(-\vp_2) K}(v_1, v_2)\I_{\xi_1+\rho u_1}(p_1) \I_{\xi_2+\rho v_1}(p_2+p_1)\,\D(u_1,u_2)\,\D(v_1, v_2)\\
   &\hspace{12cm}\times\, \gamma^{(2)}_{red}(\D p_2)\,\D p_1\\
   &=\rho \int\limits_{\R^2}\int\limits_{\R^2}\int\limits_{\R^2}\I_{O(-\vp_1)K}(u_1, u_2)\I_{O(-\vp_2) K}(v_1, v_2)\I_{\xi_1+\rho (u_1-v_1)}(p_1) \I_{\xi_2}(p_2+p_1)\,\D(u_1,u_2)\,\D(v_1, v_2)\\
    &\hspace{12cm}\times\,\gamma^{(2)}_{red}(\D p_2)\,\D p_1\,.
\end{align*}
It is easy to see that the invariance properties of the one-dimensional Hausdorff measure on $\R^2$ (also denoted by $|\cdot|_1$) yield 
$$\int\limits_{\R^1} \I_{O(-\vp_1)K}(u_1, u_2)\D u_2= |\, g(u_1,0)\cap O(-\vp_1)K \,|_1 = |\, O(\vp_1) g(u_1,0)\cap K\, |_1 =|\, g(u_1,\vp_1)\cap K\, |_1.$$
and likewise $\int_{\R^1} \I_{O(-\vp_2)K}(v_1, v_2)\,\D v_2 = |\,g(v_1,\vp_2)\cap K\,|_1$.

\medskip\noindent
Therefore, the integral $J_\rho(K,\xi_1,\vp_1,\xi_2,\vp_2)$ is equal to
\begin{align*}
    \rho & \int\limits_{\R^2}\int\limits_{\R^1}\int\limits_{\R^1} |g(u,\vp_1)\cap K|_1\, |g(v,\vp_2)\cap K|_1\, \I_{\xi_1+\rho (u - v)}(p_1)\, \I_{\xi_2}(p_2+p_1)\,\D u\,\D v\, \gamma^{(2)}_{red}(\D p_2)\,\D p_1\\
     & = \rho \int\limits_{\R^2}\int\limits_{\R^1}\int\limits_{\R^1} |g(w+v,\vp_1)\cap K|_1\, |g(v,\vp_2)\cap K|_1 \,\I_{\xi_1+\rho w}(p_1)\, \I_{\xi_2}(p_2+p_1)\,\D w\,\D v\, \gamma^{(2)}_{red}(\D p_2)\,\D p_1 \\
     & =\int\limits_{\R^2}\int\limits_{\R^1}\int\limits_{\R^1} |g(w/\rho+v,\vp_1)\cap K|_1\, |g(v,\vp_2)\cap K|_1 \,\I_{\xi_1+ w}(p_1)\, \I_{\xi_2}(p_2+p_1)\,\D w\,\D v\, \gamma^{(2)}_{red}(\D p_2)\,\D p_1 \\
 &\xrightarrow[\rho \to \infty]{} \int\limits_{\R^2}\int\limits_{\R^1}\int\limits_{\R^1} |g(v,\vp_1)\cap K|_1 \,|g(v,\vp_2)\cap K|_1\, \I_{-\xi_1+ p_1}(w)\, \I_{\xi_2-p_2}(p_1)\,\D w\D v \,\gamma^{(2)}_{red}(\D p_2)\,\D p_1\\
  &= \;|\xi_1|_1|\xi_2|_1 \gamma^{(2)}_{red}(\R^1)\int\limits_{\R^1}\I_{[\ell(\vp_1,K),r(\vp_1,K)]}(v)\, |g(v,\vp_1)\cap K|_1\, \I_{[\ell(\vp_2,K),r(\vp_2,K)]}(v)\,|g(v,\vp_2)\cap K|_1 \,\D v, 
\end{align*}
where the interval $[\ell(\vp_i,K),r(\vp_i,K)] = \{v \in \R^1: g(v, \vp_i)\cap K\neq \emptyset\}$ coincides with the orthogonal projection of $O(-\vp_i)K$ on the $v$-axis for $i=1,2$. 
To justify the above limit we have used that $|g(w/\rho+v,\vp_1)\cap K|_1 \le \mbox{diam}(K)$ so that Lebesgue's dominated convergence theorem can be applied. Furthermore, it is easily seen that
\begin{equation}\label{eq4.5}
|\,J_\rho(K,\xi_1,\vp_1,\xi_2,\vp_2)\,| \le \mbox{diam}(K)\,|K|_2\,|\xi_1|_1\,|\xi_2|_1\,\|\gamma_{red}^{(2)}\|_{TV}\,.
\end{equation}
Hence, the limit of \eqref{eq3.3}, i.e. limit of $\la\,\E J_\rho(K,\Xi_1,\Phi_1,\Xi_2,\Phi_2)$ as $\rho \to \infty\,$, exists and can be expressed by using the independence assumptions as follows:
\begin{equation*}
    \la\,(\E |\Xi_0|_1)^2 \gamma_{red}^{(2)}(\R^1)\,\int\limits_{\R^1}\big(\,\E\I_{[\ell(\Phi_0,K),r(\Phi_0,K)]}(v)\,|g(v,\Phi_0)\cap K|_1\,\big)^2\D v.
\end{equation*}
Note that the indicator function $\I_{[\ell(\Phi_0,K),r(\Phi_0,K)]}(\cdot)$ can be omitted since the range of integration w.r.t. $v$ is well-defined.
 \qed
\end{proo}

\begin{proo}[Lemma \ref{lemma6}]
With the abbreviation $\Xi_0 = [-R_0,R_0]$ we obtain that
\begin{equation*}
J_\rho(K):=\,\rho\int\limits_K \int\limits_K \int\limits_{\R^1} \vw_{\rho x,\rho y}(p)\,\D p\,\D x\,\D x=\rho \int\limits_K \int\limits_K \int\limits_{\R^1} \P\big(p \in \Xi_0 \cap (\Xi_0+\rho \langle v(\Phi_0),y-x\rangle)\big)\,\D p\, \D x \,\D y 
\end{equation*}   
    
\begin{align*} 
=\, & \rho \int\limits_{\R^2}
    \I_{K\oplus(-K)}(y)\,|K\cap(K-y)|_2\;\E|\Xi_0 \cap (\Xi_0 +\rho\,\langle v(\Phi_0),y\rangle)\,|_1\,\D y\\
 =\, & \rho \int\limits_0^{2\pi}\int\limits_0^{\infty} \I_{K\oplus(-K)}(s\, v(\psi))\,|K\cap(K-s \,v(\psi))|_2 \;\E|\Xi_0\cap (\Xi_0+\rho\, s \cos(\Phi_0-\psi))|_1\, s\, \D s\,\D \psi\,,\\
\nonumber\\ 
        &\text{where we have substituted $y=s\, v(\psi)$ with $v(\psi) = (\cos \psi, \sin \psi)^T$ and with}\\
        &\text{$r_K(\psi)={\rm max}\{s\ge 0:s\,v(\psi)\in K\oplus(-K)\}\;(\,=\,r_K(\psi\pm\pi)$ due to symmetry reasons)}
\end{align*} 

\begin{align*}        
   =\,& \rho\, \int\limits_0^{2\pi}\int\limits_0^{r_K(\psi)}|K\cap(K-s\, v(\psi))|_2 \,\E|\,\Xi_0 \cap (\Xi_0+\rho\, s\, \cos(\Phi_0-\psi))|_1\, s\, \D s\, \D \psi\\\
  =\, &\rho \,\E\int\limits_{-\Phi_0}^{2\pi-\Phi_0}\int\limits_0^{r_K(\psi+\Phi_0)}|K\cap(K-s\, v(\psi+\Phi_0))|_2 \,\E|\,\Xi_0 \cap (\Xi_0 + \rho \,s\, \cos(\psi))|_1\, s\, \D s\, \D \psi\,,\\
  &\text{where we have used the independence of $\Phi_0$ and $R_0$}
\end{align*} 
 
\begin{align*}
  =\,& \rho\, \E\int\limits_{0}^{2\pi}\int\limits_0^{r_K(\psi+\Phi_0)}|K\cap(K-s\, v(\psi+\Phi_0))|_2\, \E|\,\Xi_0 \cap (\Xi_0+\rho\, s\, \cos(\psi))|_1\, s\, \D s\, \D \psi\,,\\
&\mbox{where we have used $\int_{-\Phi_0}^0(\cdots){\rm d}\psi = \int_{2\pi-\Phi_0}^{2\pi}(\cdots){\rm d}\psi$ due to  $v(\psi) = v(\psi+2\pi)$}\\
 =\,& 2\,\rho\, \E\int\limits_0^{\pi}\int\limits_0^{r_K(\psi+\Phi_0)}|K\cap(K+s\, v(\psi+\Phi_0))|_2\, \E|\,\Xi_0 \cap (\Xi_0+\rho\, s\, \cos(\psi))|_1\, s\, \D s\, \D \psi\,,\\ 
 &\mbox{where we have used $v(\psi+\pi)= -v(\psi)$ and the shift-invariance of $|\cdot|_1$ as well as} \\
 &\mbox{the motion-invariance of $|\cdot|_2\,$.} 
\end{align*}

\noindent
By definition of $r_K(\psi)$ we have $s > r_K(\psi)$ iff $s\,v(\psi) \notin K\oplus(-K)$ iff $K \cap (K+s\,v(\psi))= \emptyset$. Thus, the inner integral $\int_0^{r_K(\psi+\Phi_0)}$ in the above double integral can be replaced by $\int_0^\infty$ showing that

\begin{align*}
J_\rho(K) &= 2\,\rho\, \E\int\limits_0^{\pi}\int\limits_0^\infty|K\cap(K+s\,v(\psi+\Phi_0))|_2\;\E|\,\Xi_0 \cap (\Xi_0+\rho\, s\, \cos(\psi))|_1\, s\, \D s\, \D \psi\\   
  &=   2\,\rho \int\limits_1^{-1}\E \int\limits_0^{\infty}|K\cap(K+s\,v(\arccos(y)+\Phi_0))|_2\;\E|\Xi_0 \cap (\Xi_0+\rho\, s\, y)|_1\,s\,\D s\, \frac{(-1)\D y}{\sqrt{1-y^2}} \\
 & \text{by substituting $y = \cos(\psi) \in [-1,1]$ so that $\psi = \arccos(y)\;,\;(\arccos(y))'= -\frac{1}{\sqrt{1-y^2}}$}
\end{align*}

\begin{align*} 
&= 2\,\rho \int\limits_{-1}^{1}\E \int\limits_0^{\infty}|K\cap(K+s\, v(\arccos(y)+\Phi_0))|_2\;\E|\Xi_0\cap(\Xi_0+\rho\, s\, y)|_1 \,\frac{s\,\D s\,\D y}{\sqrt{1-y^2}}\\
&=2 \,\E \int\limits_{0}^{\infty} \int\limits_{-s}^{s}|K \cap (K+s\, v(\arccos(\frac{z}{s})+\Phi_0))|_2\;\E|\Xi_0\cap(\Xi_0+\rho\, z)|_1 
 \,\frac{\rho\,\D z\,s\,\D s}{\sqrt{s^2-z^2}}\\
 &\text{by substituting $z = s\,y \in [-s,s]$ so that $y=z/s$ and changing the order of integration.}
\end{align*}

\bigskip\noindent
Interchanging again the integration over $z$ and $s$, we can proceed with the abbreviation
$$ h(s,z,\Phi_0):= s\,v(\arccos\big(\frac{z}{s}\big)+\Phi_0)=\big(z \cos \Phi_0-\sqrt{s^2-z^2}\sin\Phi_0, z\sin\Phi_0+\sqrt{s^2-z^2}\cos\Phi_0\big)\,,$$
where $0 \le \|h(s,z,\Phi_0)\|=s \le r_K := \max\{r_K(\vp): 0 \le \vp \le \pi\}\le \mbox{diam}(K)$, leading to
\begin{align*}
J_\rho(K) &= 2\,\E\int\limits_{\R^1}\int\limits_{|z|}^{r_K}|K\cap(K+h(s,z,\Phi_0))|_2\;\E|\Xi_0\cap(\Xi_0+\rho z)|_1\;\frac{s\,\D s\,\rho\,\D z}{\sqrt{s^2-z^2}}\\
    &=2\,\E\int\limits_{\R^1}\int\limits_{|u|/\rho}^{r_K}|K\cap(K+h(s,u/\rho,\Phi_0))|_2\;\E|\Xi_0\cap(\Xi_0+u)|_1\;\frac{s\,\D s\,\D u}{\sqrt{s^2-(u/\rho)^2}}\\
    &\text{by substituting $u = \rho\,z$ so that $z = u/\rho$}\\
    &\xrightarrow[\rho \to \infty]{}2 \,\E\int\limits_{\R^1}\int\limits_{0}^{r_K}|K\cap(K+s\, v(\Phi_0+\pi/2))|_2\;\E|\Xi_0\cap(\Xi_0+u)|_1\,\D s\,\D u.
\end{align*} 

\bigskip\noindent   
In the last line, we could apply Lebesgue's dominated convergence theorem since 
\begin{equation}\label{eq4.6}
\int\limits_{|u|/\rho}^{r_K}|K\cap(K+h(s,u/\rho,\Phi_0))|_2\;\frac{s\,\D s}{\sqrt{s^2-(u/\rho)^2}}\le |K|_2\,\frac{1}{2}\,\int\limits_0^{r_K^2-\frac{u^2}{\rho^2}}\frac{\D t}{\sqrt{t}} \le \,|K|_2\,\mbox{diam}(K)\,.
\end{equation}

\bigskip\noindent
Further, we have used the continuity of the function $z \mapsto h(s,z,\vp)$, $\,\arccos(0)=\pi/2$ and $h(s,0,\vp) = s\, v(\vp+\pi/2) = s\,(-\sin\vp,\cos\vp)^T\;(\,= -s\, v(\vp-\pi/2))$ and the relation $\int_{\R^1}|\Xi_0\cap(\Xi_0+u)|_1\,\D u = |\Xi_0|^2_1 = 4\,R^2_0$ combined with a multiple application of Fubini's theorem. Finally, we arrive at
\begin{equation*}
      J_\rho(K) = \rho \int\limits_K \int\limits_K \int\limits_{\R^1}\vw_{\rho x_1,\rho x_2}(p)\,\D p\,\D x_1\,\D x_2 
      \xrightarrow[\rho \to \infty]{}2\,\E |\Xi_0|^2_1 \int\limits_0^{\pi}\int\limits_0^{r_K(\vp \pm \frac{\pi}{2})}\big|K \cap \big(K+s v \big(\vp \pm \frac{\pi}{2}\big)\big)\big|_2 \D s \D G(\vp).
\end{equation*}\qed
\end{proo}
\bigskip

\medskip\noindent
\begin{proo}[Theorem \ref{theorem2}]

In view of Lemma \ref{Lemma2} and Definition \ref{notation} we can state the equality 
$$
\rho^{-3}\,\V\big(|\Xi \cap \rho K|_2\big) = \int\limits_K\int\limits_K \rho\,\Big(\,G_P\big[1-\uw_{\rho x_1,\rho x_2}(\cdot)\big]-G_P\big[1-w_{\rho x_1}(\cdot)\big]G_P\big[1-w_{\rho x_2}(\cdot)\big]\,\Big)\D x_1 \D x_2\,.
$$
Instead to use the factorial moment expansion of the pgf's $G_P[1-\uw_{\rho x_1, \rho x_2}], G_{P}[1-w_{\rho x_1}]$ and $G_{P}[1-w_{ \rho x_2}]$ as in \eqref{eq2.7} and \eqref{eq2.8}, we first rewrite the integrand of the right-hand side of the foregoing equality as follows: 
\begin{align}\label{eq4.7}
&\;\rho\,\Big(\,G_P\big[1-\uw_{\rho x_1,\rho x_2}(\cdot)\big]-G_P\big[1-w_{\rho x_1}(\cdot)\big]G_P\big[1-w_{\rho x_2}(\cdot)\big]\Big)
= G_P\big[1-w_{\rho x_1}(\cdot)\big]\;G_P\big[1-w_{\rho x_2}(\cdot)\big]\,\nonumber\\
& \times \,\rho\,\Big(\,\exp\Big\{\log G_P[1-\uw_{\rho x_1,\rho x_2}(\cdot)] - \log G_P[1-w_{\rho x_1}(\cdot)] - \log G_P[1-w_{\rho x_2}(\cdot)]\,\Big\} - 1 \Big)\,.
\end{align}
In order to evaluate the exponent in line \eqref{eq4.7} we use an expansion of $\log G_P\big[1-w(\cdot)\big]$ in terms of the  factorial cumulant measures $\gamma^{(k)}$ of $\Psi \sim P$, see \eqref{CumulantMeasureReprez}, which is as follows:
\begin{equation}\label{eq4.8}
\log G_P\big[1-w(\cdot)\big] = \sum_{k=1}^{\infty}\frac{(-1)^k}{k!}\int_{\R^k}\prod_{j=1}^k w(p_j)\,\gamma^{(k)}(\D(p_1,\ldots,p_k))\;,
\quad\mbox{see \cite{DVJ03},\; p.146}\,,
\end{equation}
provided the sum in \eqref{eq4.8} is convergent. In what follows we will show that  
\begin{equation}\label{eq4.9}
\overline{\lim_{\rho \to \infty}}\rho\int\limits_K\int\limits_K \Big| \log G_P[1-\uw_{\rho x_1,\rho x_2}(\cdot)] - \log G_P[1-w_{\rho x_1}(\cdot)] - \log G_P[1-w_{\rho x_2}(\cdot)]\,\Big|\,\D x_1 \D x_2 < \infty\,.
\end{equation} 

\smallskip\noindent
Before proving this, we note that the relation  \eqref{eq2.15} implies that 
\begin{equation}\label{eq4.10}
\lim_{\rho\to\infty} G_P[1-w_{\rho x}(\cdot)] =\sum_{k=0}^{m-1}\frac{(-1)^k\,(\la\,\E|\Xi_0|_1)^k}{k!} + \theta \frac{(\la\,\E|\Xi_0|_1)^m}{m!}\xrightarrow[m \to \infty]{} \exp\{-\la\,\E |\Xi_0|_1\} 
\end{equation} 
for some $\theta \in [-1,1]$ uniformly for all $x \ne \mathbf{o}$. 
Furthermore, it is rapidly seen that the limit \eqref{eq4.1} (which has been proved under the assumptions of Theorem 1)  holds if and only if 
$$
\lim_{\rho\to\infty}\Big(\,\log G_P[1-\uw_{\rho x_1,\rho x_2}(\cdot)] - \log G_P[1-w_{\rho x_1}(\cdot)] - \log G_P[1-w_{\rho x_2}(\cdot)]\,\Big) = 0
$$
for distinct points $x_1,x_2\in K\setminus\{{\bf o}\}\,$. Finally, the latter limit combined with \eqref{eq4.9} proves the equality

\begin{align*}
&\lim_{\rho\to\infty}\rho\int\limits_K\int\limits_K\Big(\exp\Big\{\log G_P[1-\uw_{\rho x_1,\rho x_2}(\cdot)]-\log G_P[1-w_{\rho x_1}(\cdot)] - \log G_P[1-w_{\rho x_2}(\cdot)]\Big\}-1\Big) \D x_1 \D x_2\nonumber\\
&= \lim_{\rho\to\infty}\rho\int\limits_K\int\limits_K\Big(\log G_P[1-\uw_{\rho x_1,\rho x_2}(\cdot)]-\log G_P[1-w_{\rho x_1}(\cdot)]-\log G_P[1-w_{\rho x_2}(\cdot)]\Big)\D x_1\D x_2.
\end{align*}
The equality of both limits results from the inequality $|e^x - 1 - x| \le \frac{x^2}{2}\,e^{\max(x,0)}$ and Lebesgue's dominated convergence theorem.

\bigskip\noindent
Combining the latter equality with the \eqref{eq4.7}, \eqref{eq4.9}, \eqref{eq4.10} and the integral representation of 
$\rho^{-3}\,\V\big(|\Xi \cap \rho K|_2\big)$ at the very beginning of the proof of Theorem 2  we can state the relation

\begin{align}\label{eq4.11}
\lim_{\rho\to\infty}\rho^{-3}\,\V\big(|\Xi \cap \rho K|_2\big)  = e^{-2 \la\,\E |\Xi_0|_1} &\,\lim_{\rho\to\infty}\int\limits_K\int\limits_K \rho\,\Big(\log G_P[1-\uw_{\rho x_1,\rho x_2}(\cdot)]\\
& - \log G_P[1-w_{\rho x_1}(\cdot)] - \log G_P[1-w_{\rho x_2}(\cdot)]\Big)\D x_1 \D x_2\,.\nonumber
\end{align}
By using the expansion \eqref{eq4.8} the double integral on the right-hand side of \eqref{eq4.11} takes the form
\begin{equation*}
\int\limits_K\int\limits_K \rho\,\Big(\log G_P[1-\uw_{\rho x,\rho y}(\cdot)] - \log G_P[1-w_{\rho x}(\cdot)] - \log G_P[1-w_{\rho y}(\cdot)]\Big)\D x \D y = \sum_{n=1}^{\infty}\frac{(-1)^nT_n^{(\rho)}(K)}{n!}\,,\\
\end{equation*}
where $T_n^{(\rho)}(K)$  for $n \in \N$ is defined by
\begin{equation}\label{eq4.12}
T_n^{(\rho)}(K):=\int\limits_K\int\limits_K\int\limits_{\R^n}\rho\,\bigg(\prod_{j=1}^n \uw_{\rho x,\rho y}(p_j)\,-\,\prod_{j=1}^n w_{\rho x}(p_j) -\prod_{j=1}^n w_{\rho y}(p_j)\bigg)\,\gamma^{(n)}(\D(p_1,\ldots,p_n))\,\D x \,\D y\,.
\end{equation}

\smallskip\noindent 
Since $\gamma^{(1)}(\D p) = \la\,\D p$ and $\uw_{\rho x,\rho y}(p)- w_{\rho x}(p)-w_{\rho y}(p) =  -\vw_{\rho x,\rho y}(p)\,$,  we get   
\begin{equation*}
 - T_1^{(\rho)}(K) = \la\,\int\limits_K\int\limits_K \int\limits_{\R^1}\rho\,\vw_{\rho x,\rho y}(p)\,\D p\,\D x \,\D y = \la\,J_\rho(K)\xrightarrow[m \to \infty]{} 2\,\la\,\E| \Xi_0 |_1^2\,C_2^{G,K}\,,
\end{equation*}
where the limit is just the assertion of Lemma \ref{lemma6}. The above proof of Lemma \ref{lemma6} reveals that  $|T_1^{(\rho)}(K)| \le \la \,J_\rho(K) \le 2\,\la \,\E| \Xi_0 |_1^2\,| K |_2\,\mbox{diam}(K)$. In the next step we derive  a uniform bound of $T_2^{(\rho)}(K)$  as well as its the limit  as $\rho \to \infty\,$. For doing this, we rewrite 
\begin{align*}
\prod_{j=1}^2\uw_{\rho x,\rho y}(p_j)-\prod_{j=1}^2 w_{\rho x}(p_j)-\prod_{j=1}^2w_{\rho y}(p_j)&= w_{\rho x}(p_1)w_{\rho y}(p_2)+w_{\rho y}(p_1)w_{\rho x}(p_2)-\uw_{\rho x,\rho y}(p_1)\vw_{\rho x,\rho y}(p_2)\\
&-\vw_{\rho x,\rho y}(p_1)\big( w_{\rho x}(p_2)+w_{\rho y}(p_2)\,\big)
\end{align*}
and by regarding the symmetry in $x,y$ and $p_1,p_2$ we get
\begin{align}\label{eq4.13}
T_2^{(\rho)}(K)& = \rho \int\limits_K\int\limits_K \int\limits_{\R^2}
\Big(\,\prod_{j=1}^2\uw_{\rho x,\rho y}(p_j)\,-\,\prod_{j=1}^2 w_{\rho x}(p_j) -\prod_{j=1}^2 w_{\rho y}(p_j)\,\Big)\gamma^{(2)}(\D(p_1,p_2))\,\D x \,\D y
\nonumber\\
&=\rho \int\limits_K\int\limits_K \int\limits_{\R^2}\Big(2\,w_{\rho x}(p_1)\,w_{\rho y}(p_2)\,-\,\big(\,\uw_{\rho x,\rho y}(p_2) + 2\,w_{\rho x}(p_2)\,\big)\,\vw_{\rho x,\rho y}(p_1)\,\Big)\gamma^{(2)}(\D(p_1,p_2))\,\D x \,\D y\nonumber\\
&= 2\,\rho \int\limits_K\int\limits_K \int\limits_{\R^2}\,w_{\rho x}(p_1)\,w_{\rho y}(p_2)\gamma^{(2)}(\D(p_1,p_2))\,\D x \,\D y + {\widetilde T_2}^{(\rho)}(K)\,,
\end{align}
where
\begin{align}\label{eq4.14}
|\,{\widetilde T_2}^{(\rho)}(K)\,| & \;\leq 3\,\la\,\rho\,\int\limits_K\int\limits_K \int\limits_{\R^1}\int\limits_{\R^1}\vw_{\rho x,\rho y}(p_1)\,w_{\rho x}(p_2+p_1)\,|\gamma^{(2)}_{red}|\,(\D p_2)\,
\D p_1\,\D x\,\D y \nonumber\\
& \;= 3\,\la\,\rho\,\int\limits_K\int\limits_K \int\limits_{\R^1}\vw_{\rho x,\rho y}(p_1)\,\E \big|\gamma^{(2)}_{red}\big|\big(\Xi_0+\rho\,\langle v(\Phi_0),x\rangle -p_1\big)\,\D p_1\,\D x\,\D y\,.
\end{align}

\smallskip\noindent
Clearly, we have $\;\infty > \big| \gamma^{(2)}_{red} \big|(\R^1)\ge  \E \big|\gamma^{(2)}_{red}\big|\big(\Xi_0+\rho\,\langle v(\Phi_0),x\rangle -p_1\big)\xrightarrow[\rho \to \infty]{} 0$ for $x \ne \bf{o}\,$. Together with the arguments used in the proof of Lemma \ref{lemma6}, among them the uniform estimate $J_\rho(K)\le 2\,\E|\Xi_0|^2\,|K|_2\,\mbox{diam}(K)$, it follows that ${\widetilde T_2}^{(\rho)}(K)\xrightarrow[\rho \to \infty]{} 0\,$. Finally, Lemma \ref{lemma5} and \eqref{eq4.13} show that

$$
\frac{T_2^{(\rho)}(K)}{2} \xrightarrow[\rho \to \infty]{} \la\,(\E\,|\,\Xi_0\,|_1)^2\,\gamma_{red}^{(2)}(\R^1)\,\int\limits_{\R^1}\,\big(\,\E\,|\,g(p,\Phi_0)\cap K\,|_1\,\big)^2\,\D p = \la\,(\E\,|\,\Xi_0\,|_1)^2\,\gamma_{red}^{(2)}(\R^1)\,C_1^{G,K}\,.
$$

\smallskip\noindent
In addition, we can derive a uniform bound of $T_2^{(\rho)}(K)$. From  \eqref{eq4.14} and the above bound of $T_1^{(\rho)}(K)$ we get that  
$|{\widetilde T_2}^{(\rho)}(K)| \le 3\,\|\gamma_{red}^{(2)}\|_{TV}\,| T_1^{(\rho)}(K)| \le  6\,\la\,|K|_2\,\mbox{diam}(K)\,\|\gamma_{red}^{(2)}\|_{TV}\,\E|\Xi_0|_1^2\,$. Hence, we see from \eqref{eq4.5} and \eqref{eq4.12} 
that, for two independent pairs $(\Xi_i,\Phi_i)\,,i=1,2,$ with the same disribution as $(\Xi_0,\Phi_0)$, the following estimate holds:
\begin{equation*}
|T_2^{(\rho)}(K)|  \le 2\,\la\,|J_{\rho}(K,\Xi_1,\Phi_1,\Xi_2,\Phi_2)| + {\widetilde T_2}^{(\rho)}(K)\, \le\, 8\,\,\la\,\,|K|_2\,\,\mbox{diam}(K)\,\,\E|\Xi_0|_1^2\,\,\|\gamma_{red}^{(2)}\|_{TV}\,.
\end{equation*}

\smallskip\noindent
Obviously, the limit \eqref{eq3.2} coincides with $\lim_{\rho \to \infty}(\,-T_1^{(\rho)}(K)+\frac{1}{2}\,T_2^{(\rho)}(K)\,)\,$. Thus, the proof of Theorem 2 is accomplished if we show that

\begin{equation}\label{eq4.15}
\lim_{\rho \to \infty} T_n^{(\rho)}(K) = 0\;\;\; \mbox{and} \;\;\; 
\sup_{\rho\ge 1} \frac{| T_n^{(\rho)}(K) |}{n!}\le C_n^K \;\;\; \mbox{for} \;\; n\ge 3 \;\;\; \mbox{such that} \;\;\; \sum_{n\ge 3}C_n^K < \infty\,. 
\end{equation}

\medskip\noindent
This means we have to find suitable upper bounds of the integrals \eqref{eq4.12} for each $n \ge 3$ which are uniform w.r.t. $\rho$ and  disappear as $\rho \to \infty\,$. Using the reduced factorial cumulant measures $\gamma_{red}^{(n)}$ defined (in differential notation) by $\gamma^{(n)}(\D(p_1,\ldots,p_n))=\la\,\gamma_{red}^{(n)}((\D p_i-p_j:i\ne j))\,\D p_j$ for any $j=1,\ldots,n\,$, the boundedness of the total variation measure $|\gamma_{red}^{(n)}|(\cdot)$ on $\R^{n-1}$ and obvious relations

\begin{align*}
& \prod_{i=1}^n\big(w_{\rho x}(p_i)+w_{\rho y}(p_i)\big)-\prod_{i=1}^n w_{\rho x}(p_i)-\prod_{i=1}^n w_{\rho y}(p_i) = \sum_{k=1}^{n-1}\sum\limits_{1\le i_1<\cdots<i_k\le n}\,\prod_{\ell=1}^k w_{\rho x}(p_{i_\ell})
\prod_{j=1\atop {j \ne i_1,\ldots,i_k}}^n w_{\rho y}(p_j)\\
& \prod_{i=1}^n\big(w_{\rho x}(p_i)+w_{\rho y}(p_i)\big)-\prod_{i=1}^n\uw_{\rho x,\rho y}(p_i) = \sum_{k=1}^n\,\vw_{\rho x,\rho y}(p_k)\,\prod_{i=1}^{k-1}\uw_{\rho x,\rho y}(p_i)\,\prod_{j=k+1}^n\big(w_{\rho x}(p_j)+w_{\rho y}(p_j)\big)\\
&\phantom{\prod_{i=1}^n\big(w_{\rho x}(p_i)+w_{\rho y}(p_i)\big)-\prod_{i=1}^n\uw_{\rho x,\rho y}(p_i)} 
\le \sum_{k=1}^n\,\vw_{\rho x,\rho y}(p_k)\,\prod_{j=1\atop {j \ne k}}^n \big(w_{\rho x}(p_j)+w_{\rho y}(p_j)\big)\,,
\end{align*}
we obtain the following estimates
 
\begin{align*}
&\bigg|\,\int\limits_K\int\limits_K \int\limits_{\R^n} \rho \Big(\,\prod_{i=1}^n\big(w_{\rho x}(p_i)+w_{\rho y}(p_i)\big) -\prod_{i=1}^n w_{\rho x}(p_i)-\prod_{i=1}^n w_{\rho y}(p_i)\,\Big)\gamma^{(n)}(\D(p_1,\ldots,p_n))\,\D x \,\D y\,\bigg|\qquad\qquad\nonumber\\
\nonumber\\
&=\bigg|\,\sum_{k=1}^{n-1}{n \choose k}\int\limits_K\int\limits_K \int\limits_{\R^n} \rho\,\prod_{i=1}^k w_{\rho x}(p_i)\,\prod_{j=k+1}^n w_{\rho y}(p_j) \,\gamma^{(n)}(\D(p_1,\ldots,p_n))\,\D x \,\D y\,\bigg|
\end{align*}

\begin{align}\label{eq4.16}
&\leq T_{n,1}^{(\rho)}(K) := \la\sum_{k=1}^{n-1}{n \choose k}\int\limits_K\int\limits_K\int\limits_{\R^1}\rho\,w_{\rho x}(p_1)\int\limits_{\R^{n-1}}\prod_{i=2}^k w_{\rho x}(p_i+p_1)\,\prod_{j=k+1}^n w_{\rho y}(p_j+p_1)\qquad\qquad\qquad\qquad\quad\nonumber\\
& \phantom{T_{n,1}^{(\rho)}(K)\,\sum_{k=1}^n{n-1 \choose k-1}\int\limits_K\int\limits_K\int\limits_{\R^1}}\qquad\qquad\qquad\qquad\qquad
 \times\,\big|\gamma_{red}^{(n)}\big|(\D(p_2,\ldots,p_n))\,\D p_1\,\D x\,\D y
\end{align}
and 
\begin{align}\label{eq4.17}
&\phantom{000} \bigg|\,\int\limits_K\int\limits_K \int\limits_{\R^n} \rho \Big(\,\prod_{i=1}^n\big(w_{\rho x}(p_i)+w_{\rho y}(p_i)\big) -\prod_{i=1}^n\uw_{\rho x,\rho y}(p_i)\,\Big)\gamma^{(n)}(\D(p_1,\ldots,p_n))\,\D x \,\D y\,\bigg|\nonumber\\
& \leq \la\,n\,\int\limits_K\int\limits_K\int\limits_{\R^1}\rho\,\vw_{\rho x,\rho y}(p_1)\int\limits_{\R^{n-1}}\,\prod_{j=2}^n\big(w_{\rho x}(p_j+p_1)+w_{\rho y}(p_j+p_1)\big)\,\big|\gamma_{red}^{(n)}\big|(\D(p_2,\ldots,p_n))\,\D p_1\,\D x \,\D y \nonumber\\
&\leq T_{n,2}^{(\rho)}(K) :=\la\,n\,\sum_{k=1}^n\,{n-1 \choose k-1}\int\limits_K\int\limits_K\int\limits_{\R^1} \rho\,\vw_{\rho x,\rho y}(p_1)\int\limits_{\R^{n-1}}\prod_{i=2}^k w_{\rho x}(p_i+p_1)\,\prod_{j=k+1}^n w_{\rho y}(p_j+p_1)\nonumber\\
&\phantom{T_{n,2}^{(\rho)}(K)\,\sum_{k=1}^n{n-1 \choose k-1}
\int\limits_K\int\limits_K\int\limits_{\R^1} \int\limits_\R^{n-1}\rho\,\qquad}\qquad\qquad\quad\times\,\big|\gamma_{red}^{(n)}\big|(\D(p_2,\ldots,p_n))\,\D p_1\,\D x\,\D y\,.
\end{align}

\smallskip\noindent
Obviously, we have $| T_n^{(\rho)}(K) | \le T_{n,1}^{(\rho)}(K) + T_{n,2}^{(\rho)}(K)$ for $n \ge 3\,$.
Let us first, rewrite the integral terms in \eqref{eq4.16}. For this purpose 
we introduce the abbreviaton 
$$
I_{n,k}^{(\rho)}(K) := \int\limits_K\int\limits_K\int\limits_{\R^1}\rho\,w_{\rho x}(p_1)\int\limits_{\R^{n-1}}\prod_{i=2}^k w_{\rho x}(p_i+p_1)\,\prod_{j=k+1}^n w_{\rho y}(p_j+p_1)\,\big|\gamma_{red}^{(n)}\big|(\D(p_2,\ldots,p_n))\,\D p_1\,\D x\,\D y
$$
for $k=2,\ldots,n-1\,$.\\
As in \eqref{eq4.4} we substitute $x = O(\Phi_1)\,u$ and $y = O(\Phi_n)\,w$ with $O(\cdot)$ as defined in \eqref{orthmatrix}. Since $O^{-1}(\vp)= O(-\vp)$ and $\det(O(\vp))=1$ it follows that $u = O(-\Phi_1)\,x\,$,  $w = O(-\Phi_n)\,y$    and  $\langle v(\Phi_i),x\rangle = \langle v(\Phi_i),O(\Phi_1)\,u \rangle =  \langle v(\Phi_i-\Phi_1),u\rangle$ for $i=1,\ldots,k$ and $\langle v(\Phi_j),y\rangle = \langle v(\Phi_j-\Phi_n),w \rangle$ for $j=k+1,\ldots,n\,$.  Note that $\langle v(\Phi_1),x\rangle =u_1$ and $\langle v(\Phi_n),y\rangle =w_1$ for $u=(u_1,u_2)^T$ and $w=(w_1,w_2)^T$, respectively.

\bigskip\noindent
Similarly as in the proof of Lemma 3  we introduce  independent copies     
$(R_1,\Phi_1),\ldots,(R_n,\Phi_n)$  of the random vector $(R_0,\Phi_0)$ and 
independent copies $\Xi_1,\ldots,\Xi_n$ of the random intervall $\Xi_0=[-R_0,R_0]\,$. Then the product $w_{\rho x}(p)\prod_{i=2}^k w_{\rho x}(p_i+p_1)\,\prod_{j=k+1}^n w_{\rho y}(p_j+p_1)$ can be expressed as expectation
\begin{equation*}
\E\Big(\I_{\Xi_1+\rho\langle\Phi_1,x\rangle}(p_1)\prod_{i=2}^k \I_{\Xi_i+\rho\langle v(\Phi_i),x\rangle}(p_i+p_1)\prod_{j=k+1}^n \I_{\Xi_j+\rho\langle v(\Phi_j),y\rangle}(p_j+p_1)\Big)
\end{equation*}
which together with the above transformations of $x, y \in \R^2$ and Fubini's theorem allows us to write $I_{n,k}^{(\rho)}(K)$ in the form 

\medskip
\begin{align*}
\phantom{=}&\;\E\int\limits_{\R^2}\int\limits_{\R^2}\int\limits_{\R^1}\int\limits_{\R^{n-1}}\rho\,\Big(\I_{\Xi_1+\rho\,u_1}(p_1)\,\prod_{i=2}^k \I_{\Xi_i+\rho\langle v(\Phi_i-\Phi_1),u\rangle}(p_i+p_1)\prod_{j=k+1}^{n-1} \I_{\Xi_j+\rho\langle v(\Phi_j-\Phi_n),w\rangle}(p_j+p_1)\nonumber\\
&\phantom{\quad}\times\,\I_{\Xi_n}(p_n+p_1-\rho\,w_1)\Big)\,\big|\gamma_{red}^{(n)}\big|(\D(p_2,\ldots,p_n))\,\D p_1\,\I_{O(-\Phi_1)K}(u)\,\I_{O(-\Phi_n)K}(w)\,\D (u_1,u_2)\,\D (w_1,w_2)\nonumber\\
= &\;\E\int\limits_{\R^2}\int\limits_{\R^2}\int\limits_{\R^1}\int\limits_{\R^{n-1}}\rho\,\Big(\I_{\Xi_1+\rho\,(u_1-w_1)}(p_1-\rho\,w_1)\prod_{i=2}^k \I_{\Xi_i+\rho\langle v(\Phi_i-\Phi_1),u\rangle-\rho w_1}(p_i+p_1-\rho w_1)\nonumber\\
\times &\prod_{j=k+1}^{n-1} \I_{\Xi_j+\rho\langle v(\Phi_j-\Phi_n),w\rangle-\rho w_1}(p_j+p_1-\rho w_1)\,\I_{\Xi_n}(p_n+p_1-\rho\,w_1)\Big)\,\big|\gamma_{red}^{(n)}\big|(\D(p_2,\ldots,p_n))\,\D p_1\,
\nonumber\\
\times &\,\I_{O(-\Phi_1)K}(u)\,\I_{O(-\Phi_n)K}(w)\,\D (u_1,u_2)\,\D (w_1,w_2)\nonumber\\
= &\;\E\int\limits_{\R^2}\int\limits_{\R^2}\int\limits_{\R^1}\int\limits_{\R^{n-1}}\rho\,\Big(\prod_{i=2}^k \I_{\Xi_i+\rho\langle v(\Phi_i-\Phi_1),u\rangle-\rho w_1}(p_i+p_1)\prod_{j=k+1}^{n-1} \I_{\Xi_j+\rho\langle v(\Phi_j-\Phi_n),w\rangle-\rho w_1}(p_j+p_1)\,\nonumber\\
\times &\,\I_{\Xi_1+\rho\,(u_1-w_1)}(p_1)\,\I_{\Xi_n}(p_n+p_1)\Big)\,\big|\gamma_{red}^{(n)}\big|(\D(p_2,\ldots,p_n))\,\D p_1\,\I_{O(-\Phi_1)K}(u)\,\I_{O(-\Phi_n)K}(w)\,\D (u_1,u_2)\,\D (w_1,w_2)\nonumber\\
= &\;\E\int\limits_{\R^2}\int\limits_{\R^2}\int\limits_{\R^1}\int\limits_{\R^{n-1}}\rho\,\Big(\prod_{i=2}^k \I_{\Xi_i+\rho\langle v(\Phi_i-\Phi_1),(z_1+w_1,z_2)\rangle-\rho w_1}(p_i+p_1)\prod_{j=k+1}^{n-1} \I_{\Xi_j+\rho\langle v(\Phi_j-\Phi_n),w\rangle-\rho w_1}(p_j+p_1)\nonumber\\
\times &\,\I_{\Xi_1+\rho\,z_1}(p)\,\I_{\Xi_n}(p_n+p)\Big)\,\big|\gamma_{red}^{(n)}\big|(\D(p_2,\ldots,p_n))\,\D p\,\I_{O(-\Phi_1)K}((z_1+w_1,z_2))\,\I_{O(-\Phi_n)K}(w)\,\D (z_1,z_2)\,\D (w_1,w_2)\nonumber\\
\end{align*}
\begin{align}\label{eq4.18} 
= & \E\int\limits_{\R^2}\int\limits_{\R^2}\int\limits_{\R^1}\int\limits_{\R^{n-1}}\Big(\prod_{i=2}^k \I_{\Xi_i+\langle v(\Phi_i-\Phi_1),(z_1+\rho w_1,\rho z_2)\rangle-\rho w_1}(p_i+p_1)\prod_{j=k+1}^{n-1} \I_{\Xi_j+\rho\langle v(\Phi_j-\Phi_n),w\rangle-\rho w_1}(p_j+p_1)\\
&\phantom{\qquad\qquad\qquad\qquad} \times \I_{\Xi_n}(p_n+p_1)\Big)\,\big|\gamma_{red}^{(n)}\big|(\D(p_2,\ldots,p_n))\,\I_{\Xi_1+z_1}(p_1)\,\D p_1\,\nonumber\\
&\phantom{\qquad\qquad\qquad\qquad} \times \I_{O(-\Phi_1)K}((\frac{z_1}{\rho}+w_1,z_2))\,\I_{O(-\Phi_n)K}((w_1,w_2))\,\D (z_1,z_2)\,\D (w_1,w_2)\,.
\nonumber
\end{align}

\bigskip\noindent
Replacing the two products of indicator functions in \eqref{eq4.18} by 1 leads to the following bound of $I_{n,k}^{(\rho)}(K)$ provided that $\big|\gamma_{red}^{(n)}\big|(\R^{n-1}) < \infty\,$:

\begin{align}
I_{n,k}^{(\rho)}(K) &\leq \E\int\limits_{\R^2}\int\limits_{\R^2}\int\limits_{\R^1}\int\limits_{\R^{n-1}}\Big(\I_{-\Xi_1+p_1}(z_1)\,\I_{\Xi_n-p_n}(p_1)\Big)\,\big|\gamma_{red}^{(n)}\big|(\D(p_2,\ldots,p_n))\,\D p_1\nonumber\\
&\phantom{\qquad\qquad\qquad\qquad}\times\,\I_{O(-\Phi_1)K}((\frac{z_1}{\rho}+w_1,z_2))\,\D (z_1,z_2)\,\I_{O(-\Phi_n)K}((w_1,w_2))\,\D (w_1,w_2)\nonumber
\end{align}

\begin{align}\label{eq4.19}
& = \E\int\limits_{\R^1}\int\limits_{\R^1}\int\limits_{\R^1}\int\limits_{\R^{n-1}}\Big(\I_{-\Xi_1+p_1}(z_1)\,\I_{\Xi_n-p_n}(p_1)\Big)\,\big|\gamma_{red}^{(n)}\big|(\D(p_2,\ldots,p_n))\,\D p_1 \;\nonumber\\
& \phantom{\qquad\qquad\qquad\qquad}\times\,|g(\frac{z_1}{\rho}+w_1,\Phi_1)\cap K|_1\;|g(w_1,\Phi_n)\cap K|_1 \D z_1\D w_1 \nonumber\\
& \leq \mbox{diam}(K)\,\E\int\limits_{\R^1}|g(w_1,\Phi_n)\cap K|_1\D w_1 
\int\limits_{\R^1}\int\limits_{\R^1}\I_{-\Xi_1+p_1}(z_1)\,\I_{\Xi_n-p_n}(p_1)\D z_1 \D p_1\;\big|\gamma_{red}^{(n)}\big|(\R^{n-1})\nonumber\\
\nonumber\\
& = \;\mbox{diam}(K)\;|K|_2\;\E|\Xi_1|_1\;\E|\Xi_n|_1\;\big|\gamma_{red}^{(n)}\big|(\R^{n-1})=\mbox{diam}(K)\;|K|_2\;(\E|\Xi_0|_1)^2\,\|\gamma_{red}^{(n)}\|_{TV}\,.
\end{align}

\medskip\noindent
Here, we have used arguments which have already been applied to prove \eqref{eq4.5}. On the other hand, the product of the indicator functions in the first line of \eqref{eq4.18} disappears as $\rho \to \infty\,$  $\P$-almost surely and for almost all  $(w_1,w_2),(z_1,z_2),p_1,(p_2,\ldots,p_n)\in \R^{n+4}$ w.r.t. the corresponding product measure. Therefore, again by Lebesgue's dominated convergence theorem,
\begin{equation}\label{eq4.20} 
\lim_{\rho \to \infty}I_{n,k}^{(\rho)}(K) = 0\quad\mbox{for}\quad k=2,\ldots,n\;,\;n \ge 3\,. 
\end{equation}

\smallskip\noindent
Next, we  derive a  further bound of $I_{n,k}^{(\rho)}(K)$ that depends more on the mean thickness $\E|\Xi_0|_1$ of the typical cylinder. For this, we need the Radon-Nikodym density $|\,c_{red}^{(n)}(p_2,\ldots,p_n)\,|$  of   $|\,\gamma_{red}^{(n)}(\cdot)\,|$ w.r.t.
to Lebesgue measure on $\R^{n-1}$. Hence, by using Fubini's theorem, we replace the integral \eqref{eq4.19} over $\R^{n-1}$ by two iterated integrals. The first integral over $(p_2,\ldots,p_{n-1})\in\R^{n-2}$ can be estimated by H\"older's inequality as follows :

\smallskip\noindent
\begin{align}\label{eq4.21}
& \qquad\qquad\int\limits_{\R^{n-2}}\prod_{i=2}^k \I_{\Xi_i+\langle v(\Phi_i-\Phi_1),(z_1+\rho w_1,\rho z_2)\rangle-\rho w_1-p_1}(p_i)\,\prod_{j=k+1}^{n-1} \I_{\Xi_j+\rho \langle v(\Phi_j-\Phi_n),w\rangle-\rho w_1-p_1}(p_j)
\nonumber\\
& \qquad\qquad\;\times\;\big|c_{red}^{(n)}(p_2,\ldots,p_{n-1},p_n)\big|\,\D(p_2,\ldots,p_{n-1}) \leq \nonumber\\
& \nonumber\\
& \Big(\int\limits_{\R^{n-2}}\prod_{i=2}^k \I_{\Xi_i+\langle v(\Phi_i-\Phi_1),(z_1+\rho w_1,\rho z_2)\rangle-\rho w_1-p_1}(p_i)\prod_{j=k+1}^{n-1} \I_{\Xi_j+\rho\langle v(\Phi_j-\Phi_n),w\rangle-\rho w_1-p_1}(p_j)\D(p_2,\ldots,p_{n-1})\Big)^\frac{q-1}{q}\nonumber\\
& \;\times \Big(\int\limits_{\R^{n-1}}\big|c_{red}^{(n)}(p_2,\ldots,p_{n-1},p_n)\big|^q\,\D(p_2,\ldots,p_{n-1})\,\Big)^\frac{1}{q} =  \Big(\prod_{i=2}^{n-1}|\Xi_i|_1\Big)^\frac{q-1}{q}\;\|c_{red}^{(n)}(\cdot,p_n)\|_q
\end{align}
for any $ q > 1\,$, where $\|c_{red}^{(n)}(\cdot,p_n)\|_q$ coincides with the term in front of the equal sign in \eqref{eq4.21}.  Combining the estimates \eqref{eq4.16} and \eqref{eq4.21} with $|g(p,\vp)\cap K|_1 \le \mbox{diam}(K)$ for $(p,\vp)\in \R^1\times [0,\pi]$,
$\int_{\R^1}|g(p,\vp)\cap K|_1 \D p = |K|_2$, switching the order of integration and finally applying Lyapunov's inequality we arrive at

\begin{align*}
I_{n,k}^{(\rho)}(K) &\leq \E\int\limits_{\R^1}\int\limits_{\R^1}\int\limits_{\R^1}\int\limits_{\R^1}\Big(\prod_{i=2}^{n-1}|\Xi_i|_1\Big)^\frac{q-1}{q}\;\|c_{red}^{(n)}(\cdot,p_n)\|_q\;
\I_{-\Xi_1+p_1}(z_1)\,\I_{\Xi_n-p_n}(p_1)\,\D p_1\,\D p_n\nonumber\\
& \phantom{\qquad\qquad\qquad\qquad\qquad\qquad}\times\,|g(\frac{z_1}{\rho}+w_1,\Phi_1)\cap K|_1\;|g(w_1,\Phi_n)\cap K|_1 \D z_1\D w_1 \nonumber\\
\nonumber\\
& \leq \mbox{diam}(K)\,|K|_2\, 
\E\Big(\prod_{i=2}^{n-1}|\Xi_i|_1\Big)^\frac{q-1}{q}\,\int\limits_{\R^1}\int\limits_{\R^1}\int\limits_{\R^1}\|c_{red}^{(n)}(\cdot,p_n)\|_q\,
\I_{-\Xi_1+p_1}(z_1)\,\I_{\Xi_n-p_n}(p_1)\D z_1\D p_1\D p_n\;\nonumber\\
\nonumber\\
& = \mbox{diam}(K)\,|K|_2\,\int\limits_{\R^1}\|c_{red}^{(n)}(\cdot,p)\|_q\, \D p \;\big(\,\E|\Xi_0|_1\,\big)^{\frac{n(q-1)}{q}+\frac{2}{q}}\,.
\end{align*}
Applying the same arguments as above, the estimate \eqref{eq4.20} reveals  that \eqref{eq4.20}
remains true if, instead of $\|\gamma_{red}^{(n)}\|_{TV} < \infty $, the $L^*_q-$norm  $\|c_{red}^{(n)}\|^*_q := \int_{\R^1}\|c_{red}^{(n)}(\cdot,p)\|_q\,\D p$ is finite for some $q > 1$ and $n \ge 3\,$. Hence, we have
\begin{align*}
T_{n,1}^{(\rho)}(K) = \la\sum_{k=1}^{n-1}{n \choose k}\,
I_{n,k}^{(\rho)}(K) \le \la\,\mbox{diam}(K)\,|K|_2\;(2^n-2)\;\big(\,\E|\Xi_0|_1\,\big)^{\frac{n(q-1)}{q}+\frac{2}{q}}\;\|c_{red}^{(n)}\|_q^*\,.
\end{align*}
Together with the strong  $L^*_q$-Brillinger mixing condition with $b_q^*\,(\E|\Xi_0|_1)^{1-\frac{1}{q}} < 1/2$ we get 
\begin{equation*}
\sum_{n\ge 3} \frac{T_{n,1}^{(\rho)}(K)}{n!}\le  \la\,a_q^*\,(\E|\Xi_0|_1)^{\frac{2}{q}}\,\mbox{diam}(K)\,|K|_2\sum_{n \ge 3}\big(2\,b_q^*\,(\E|\Xi_0|_1)^\frac{(q-1)}{q}\big)^n \le \frac{\la\,a_q^*\,(\E|\Xi_0|_1)^{\frac{2}{q}}\,\mbox{diam}(K)\,|K|_2}{1-2\,b_q^*\,(\E|\Xi_0|_1)^{1-\frac{1}{q}}}\,.
\end{equation*}

\medskip\noindent
Next, we derive two different bounds for the sum $T_{n,2}^{(\rho)}(K)$ defined in \eqref{eq4.17}. For this purpose, in analogy to $I_{n,k}^{(\rho)}(K)$, we need uniform bounds of 
\begin{equation*}
J_{n,k}^{(\rho)}(p) := \int\limits_{\R^{n-1}}\prod_{i=2}^k w_{\rho x}(p_i+p_1)\,\prod_{j=k+1}^n w_{\rho y}(p_j+p)\,\big|\gamma_{red}^{(n)}\big|(\D(p_2,\ldots,p_n))\,.
\end{equation*}
It is easily seen that
\begin{align*}
& J_{n,k}^{(\rho)}(p) = \E\int\limits_{\R^{n-1}}\prod_{i=2}^k \I_{\Xi_i+\rho\langle v(\Phi_i),x\rangle-p_1}(p_i)\prod_{j=k+1}^n \I_{\Xi_j+\rho\langle v(\Phi_j),y\rangle-p_1}(p_j)\,\big|\gamma_{red}^{(n)}\big|\D(p_2,\ldots,p_{n-1}) \le \big|\gamma_{red}^{(n)}\big|(\R^{n-1}) \nonumber\\
& \mbox{and, for any}\;\; q > 1 \;\mbox{such that}\;\; \|c_{red}^{(n)}\|_q < \infty,\nonumber\\
& J_{n,k}^{(\rho)}(p) = \E\int\limits_{\R^{n-1}}\prod_{i=2}^k \I_{\Xi_i+
\rho\langle v(\Phi_i),x\rangle-p}(p_i)\prod_{j=k+1}^n \I_{\Xi_j+\rho\langle v(\Phi_j),y\rangle-p}(p_j)\,c_{red}^{(n)}(p_2,\ldots,p_{n-1})\D(p_2,\ldots,p_{n-1}) \nonumber\\
&\le\E  \prod_{i=2}^n\big|\Xi_i\big|^{\frac{q-1}{q}}\Big(\int\limits_{\R^{n-1}}\big|c_{red}^{(n)}(p_2,\ldots,p_{n-1})\big|^q\D(p_2,\ldots,p_{n-1})\Big)^\frac{1}{q} \le (\E|\Xi_0|)^{(n-1)\frac{q-1}{q}}\;\|c_{red}^{(n)}\|_q\,.
\end{align*}

\smallskip\noindent
The foregoing estimates show that 
\begin{equation}\label{eq4.22} 
\lim_{\rho \to \infty}J_{n,k}^{(\rho)}(p) = 0\quad\mbox{for}\quad k=2,\ldots,n\;,\;n \ge 3\,. 
\end{equation}

\smallskip\noindent
Further, from the definition of $T_{n,2}^{(\rho)}(K)$, see  \eqref{eq4.17}, and the integral $J_\rho(K)$ introduced and estimated in the proof of Lemma 6    with the uniform upper bound $2\,\mbox{diam}(K)\,|K|_2\,\E|\Xi_0|_1^2$, 
we see that   
\begin{align*}
T_{n,2}^{(\rho)}(K) &\le \la\,n\,\sum_{k=1}^n\,{n-1 \choose k-1}\int\limits_K\int\limits_K\int\limits_{\R^1} \rho\,\vw_{\rho x,\rho y}(p_1)\,\D p_1\,\D x\,\D y\,\max_{2\le k\le n}\sup_{p \in \R^1} J_{n,k}^{(\rho)}(p)\nonumber\\
&= \la\;n\,2^{n-1}\;J_\rho(K)\,\max_{2\le k\le n}\sup_{p \in \R^1} J_{n,k}^{(\rho)}(p) \le \la\;n\,2^n\,\mbox{diam}(K)\,|K|_2\,\E|\Xi_0|_1^2\,\max_{2\le k\le n}\sup_{p \in \R^1} J_{n,k}
\end{align*}

\smallskip\noindent
Under the assumption that $\Psi \sim P$ is either strongly Brillinger-mixing with $b < 1/2$ or  strongly $L_q$-Brillinger-mixing with  $b_q\,(\E|\Xi_0|_1)^{1-\frac{1}{q}} < 1/2$ we obtain the inequalities   
\begin{align*}
&\sum_{n\ge 3}\frac{T_{n,2}^{(\rho)}(K)}{n!}\le 2\,\la\,a\,b\,\E|\Xi_0|_1^2\,\mbox{diam}(K)\,|K|_2\sum_{n \ge 3}n\,\big(2\,b\big)^{n-1} 
\le \frac{2\la\,a\,b\,\E|\Xi_0|_1^2\,\mbox{diam}(K)\,|K|_2}{(1-2\,b)^2}
\nonumber\\
\mbox{and}\qquad &\nonumber\\
&\sum_{n\ge 3}\frac{T_{n,2}^{(\rho)}(K)}{n!} \le 2\,\la\,a_q\,b_q\,\E|\Xi_0|_1^2\,\mbox{diam}(K)\,|K|_2\sum_{n \ge 3}n\big(2\,b_q\,(\E|\Xi_0|_1)^\frac{q-1}{q}\big)^{n-1}\nonumber\\
&\phantom{\sum_{n\ge 3}\frac{T_{n,2}^{(\rho)}(K)}{n!}} \le\frac{2\,\la\,a_q\,b_q\,\E|\Xi_0|_1^2\,\mbox{diam}(K)\,|K|_2}{\big( 1-2\,b_q\,(\E|\Xi_0|_1)^{1-\frac{1}{q}}\big)^2}.
\end{align*}

\bigskip\noindent
Finally, summarizing the above-proved relations \eqref{eq4.20}, \eqref{eq4.22} and the convergence  of the series $\sum_{n\ge 3}T_{n,i}^{(\rho)}(K)/n!$ for $i=1,2$  shows the validity of \eqref{eq4.15} which in turn implies \eqref{eq4.9}. Thus, the proof of Theorem 2 is complete. \qed
\end{proo}

\bigskip\noindent
\begin{rema}\label{remark6}  
Strong Brillinger-mixing with $b < 1/2$ is a rather restrictive condition for the one-dimensional point process $\Psi\sim P$. Equivalently formulated, the power series $\sum_{n=2}^\infty \frac{z^n}{n!}\,|\gamma^{(n)}(\R^{n-1})|$ is analytic in the interior of the disk $b({\bf{o}},2)$ in the complex plane. For example, the condition has been used for statistical analysis of point processes in \cite{Dav77}. The Gauss-Poisson process, see \cite{DVJ03}, Poisson cluster processes with a finite number of non-vanishing cumulant measures, see \cite{AmmTha77,AmmTha78}, and certain Neyman-Scott processes satisfy this condition. In the case that $\Psi\sim P$ is strongly $L_q-$ resp. strongly $L^*_q-$Brillinger-mixing for some $q > 1$ with $b_q > 0$ resp. $b^*_q > 0$ we can choose $\E|\Xi_0|$ sufficiently small to fulfill the assumptions of Theorem 2 which greatly expands its applicability.  
\end{rema} 

\medskip\noindent
In a separate paper we will study the asymptotic normality of 
$\rho^{-3/2}\;\big(\,|\Xi \cap \rho K|_2 - \E|\Xi \cap \rho K|_2\,\big)$ as $\rho \to \infty\,$. To achieve this goal  we have to find the conditions which allow to verify that

\begin{equation}\label{eq4.23}
\rho^{-3\,k/2}\;\C_k(|\Xi \cap \rho\,K|_2) \xrightarrow[\rho \to \infty]{} 0\;\;\mbox{for any}\;\;k \ge 3\,,
\end{equation}
where with the notation and the formulas of Chapter 2 we can use the following representation of the $k$th-order cumulant $\C_k(|\Xi \cap \rho\,K|_2)$:
\begin{equation*}
(-1)^k\,\sum_{\ell=1}^k (-1)^{\ell-1}(\ell-1)! \sum
\limits_{\substack{K_1\cup\cdots\cup K_\ell\\=\{1,\dots,k\}}}\int\limits_{(\rho\,K)^k}\,\prod_{j=1}^\ell G_P\bigl[1-{\bf P}\bigl((\cdot)\in \bigcup_{i\in K_j}(\Xi_0 + \langle v(\Phi_0),x_i\rangle)\bigr)\bigr]\,\D(x_1,\ldots,x_k)\,.
\end{equation*}
From the latter formula it is easily seen that \eqref{eq4.23} is equivalent to
\begin{align*}
\rho^{k/2}\;\sum_{\ell=1}^k \frac{(-1)^{\ell-1}}{\ell} \sum\limits_{\substack{k_1+\cdots+k_\ell=k\\k_i\geq 1, i=1\dots,\ell}}\frac{k!}{k_1!\cdots k_\ell!}\prod_{j=1}^\ell\int\limits_{K^{k_j}}\prod_{j=1}^\ell G_P\bigl[1-\uw_{\rho x_1,\ldots,\rho x_{k_j}}(\cdot)\bigr]\,\D(x_1,\ldots,x_{k_j})\xrightarrow[\rho \to \infty]{} 0
\end{align*}
for any $k \ge 3\,$.  A modification of a recursive proving technique developed in Chapter 2 of \cite{HeiSpi09} to treat the analogous problem for Poisson cylinder processes could be useful.

\section*{Acknowledgements}
\noindent
The research of Flimmel was supported by the Czech Science Foundation, project 17-00393J, by Charles University Grant Agency, project No. 466119 and by a scholarship from the German Academic Exchange Service (DAAD). She is especially thankful to the University of Augsburg for its kind hospitality and support.

\end{document}